\documentstyle[11pt]{article}
\newdimen\AAdi%
\newbox\AAbo%
\def\AAk#1#2{\setbox\AAbo=\hbox{#2}\AAdi=\wd\AAbo\kern#1\AAdi{}}%
\def\AAr#1#2#3{\setbox\AAbo=\hbox{#2}\AAdi=\ht\AAbo\raise#1\AAdi\hbox{
#3}}%
 \font\tenmsb=msbm10 at 12pt
\font\sevenmsb=msbm7 at 8pt \font\fivemsb=msbm5 at 6pt
\newfam\msbfam
\textfont\msbfam=\tenmsb \scriptfont\msbfam=\sevenmsb
\scriptscriptfont\msbfam=\fivemsb
\def\Bbb#1{{\tenmsb\fam\msbfam#1}}

\def\R{\Bbb R}

\def\B{\Bbb B}
\def\N{\Bbb N}

\title{Bubble towers for supercritical semilinear elliptic equations}

\author{Yuxin Ge,$^a$ Ruihua Jing $^{a,b}$ and Frank Pacard $^a$
\thanks {{\it E-mail addresses:} ge@univ-paris12.fr (Y. Ge), jing@univ-paris12.fr (R. Jing), pacard@univ-paris12.fr (F. Pacard)}}
\date{}
\begin{document}

\maketitle

\newtheorem{th}{Theorem}
\newtheorem{lem}{Lemma}
\newtheorem{cor}{Corollary}
\newtheorem{rem}{Remark}
\newtheorem{pro}{Proposition}
\newtheorem{defi}{Definition}
\newcommand{\beq}{\begin{equation}}
\newcommand{\beqr}{\begin{eqnarray}}
\newcommand{\noi}{\noindent}
\newcommand{\eeq}{\end{equation}}
\newcommand{\eeqr}{\end{eqnarray}}
\newcommand{\mint}{-\!\!\!\!\!\!\int}
\addtocounter{rem}{-1}

\def \bx{\hspace{2.5mm}\rule{2.5mm}{2.5mm}} \def \vs{\vspace*{0.2cm}}
\def
\hs{\hspace*{0.6cm}}
\def \ds{\displaystyle}
\def \p{\partial}
\def \Wp{W^{1, p}}
\def \e{\varepsilon}
\def \o{\omega}
\def \t{\theta}
\def \na{\nabla}
\def \O{\Omega}
\def \Bar{\B_d(x_0, \alpha R)}
\def \La{\Delta}
\def \l{\lambda}
\def \s{\sigma}
\def \ut{u^\tau}
\def\cqfd{%
\mbox{ }%
\nolinebreak%
\hfill%
\rule{2mm} {2mm}%
\medbreak%
\par%
}
\def \pr {\noindent {\it Proof:} }
\def \rmk {\noindent {\it Remark} }
\def \esp {\hspace{4mm}}
\def \dsp {\hspace{2mm}}
\def \ssp {\hspace{1mm}}

\begin{center}
$^a$ Laboratoire  d'Analyse et de Math\'ematiques Appliqu\'ees, 
CNRS UMR 8050, 
D\'epartement de Math\'ematiques, 
Universit\'e Paris XII-Val de Marne, \\61 avenue du G\'en\'eral de Gaulle, 
94010 Cr\'eteil Cedex, France\\[0cm]
$^b$ Department of Mathematics, 
East China Normal University, 
200062 Shanghai, \\P.R. of China\\
\end{center}

{\bf Abstract :} We construct positive solutions of the semilinear
elliptic problem $\Delta u+ \lambda \, u + u^p = 0$ with Dirichet
boundary conditions, in a bounded smooth domain $\Omega \subset
\R^N$ $(N\geq 4)$, when the exponent $p$ is supercritical and
close enough to $\frac{N+2}{N-2}$ and the parameter $\lambda\in\R$ is
small enough. As $p\rightarrow \frac{N+2}{N-2}$, the solutions
have multiple blow up at finitely many points which are the
critical points of a function whose definition involves Green's
function. Our result extends the result of Del Pino, Dolbeault and
Musso \cite{DDM} when $\Omega$ is a ball and the solutions are
radially symmetric.



\section{Introduction}

In this paper we consider the semilinear elliptic problem \beq
\label{eq1.1} \left\{\begin{array}{rllllll} \Delta u + \lambda \,
u  +  u^p & = & 0 \quad & \mbox{in} \quad \O
\\[3mm]
u & > &  0 \quad &\mbox{in} \quad \O \\[3mm]
u & = & 0 \quad &\mbox{on} \quad \partial\O
\end{array}\right.
\eeq where $\O$ is a bounded regular domain in $\R^N$, $N\geq 4$, the
parameter $\lambda \in\R$ and the exponent $p$ is larger than
\[
p_N:=\frac{N+2}{N-2},
\]
the critical Sobolev exponent. 

\medskip

When $p=p_N$, Brezis and Nirenberg  \cite{BN} have proved that (\ref{eq1.1}) admits a solution  provided $0<\lambda$ is less than the first eigenvalue of the Laplacian on $\O$ with $0$ Dirichlet boundary condition. Direct application of Pohozaev's identity \cite{P} shows that solutions of (\ref{eq1.1}) do not exist when $\lambda \leq 0$,  $p\geq p_N$ and $\O$ is a star-shaped domain. 

\medskip

In this paper, we are interested in the existence of solutions of (\ref{eq1.1}) in the case where $p$ is larger than the critical Sobolev exponent. When $\O$ is the unit ball it is easy to check that there exist radially symmetric positive solutions of
\[
\Delta u  +  u^p = 0 ,
\]
which have multiple blow up at the origin as the exponent $p$
tends to $p_N$ (we do not assume Dirichlet boundary condition
here). We discuss this result in section 4.  In a recent paper
\cite{DDM}, Del Pino, Dolbeault and Musso have proved that a
similar result was also true for (\ref{eq1.1}). These solutions
which have multiple blow up at some points in $\O$ will be
referred to as "bubble tree solutions". We are interested in the
existence of these bubble tree solutions when $\Omega$ is
arbitrary.

\medskip

\section{Statement of the result}

Let $G$ denote Green's function for the Laplace operator with
Dirichlet boundary condition on $\O$ and let $H$ denote Robin's
function, i.e. the regular part of Green's function. Namely
\[
G(y,z) : =  |y-z|^{2-N} - H(y,z) ,
\]
for $(y,z)\in\Omega\times\Omega$. Observe that $ \Delta_y H  =  0$
in $\Omega\times\Omega$ and $G =0$ on $\partial (\O\times \O)$.

\medskip

Given $m\in \N^*$ and ${\bf x} := (x_1,\ldots ,x_m)\in \Omega^m$,
we define the $m \times m$ matrix
\[
M ({\bf x}) : = (m_{ij})_{1\leq i, j\leq m},
\]
whose entries are given by
\begin{equation}
m_{ii} : = H(x_i,x_i) >  0 \qquad \mbox{and} \qquad m_{ij} : = -
G(x_i,x_j) < 0 , \label{eq1.2}
\end{equation}
if $i \neq j$. Let  $\rho({\bf x})$ be the least eigenvalue of
$M({\bf x})$. We agree that $\rho({\bf x}) = - \infty$, if
$x_i=x_j$ for some $i\neq j$. Finally, we define $r({\bf x})$ to
be the unique eigenvector associated to $\rho ({\bf x})$ whose
coordinates are all positive and which is normalized so that its
norm is equal to $1$ (given the signs of the entries of $M({\bf
x})$, it is easy to check that one can choose the eigenvector
corresponding to the least eigenvalues to have coordinates greater
than $0$).

\medskip

We define the open set
\[
P^+ : = \{ (\Lambda, {\bf x}) \in (\R^*_+)^m \times \Omega^m \, \}.
\]
Given $\mu \in {\R}$ and  $(\ell_1,\ldots ,\ell_m) \in \N^m$, we
define
\begin{equation}
{\cal F}_\mu : (\R^*_+)^m \times \Omega^m \longrightarrow \R ,
\label{eq1.3}
\end{equation}
by
\[
{\cal F}_\mu (\Lambda, {\bf x} ) : = \Lambda \, M({\bf x})
 \, {}^t\Lambda - \mu \, C_N^{(1)} \, \sum_{i=1}^m \Lambda_i^{4\over
N-2}+ C_N^{(2)}\sum_{i=1}^m \ell_i \, \log \Lambda_i ,
\]
where $\Lambda =  (\Lambda_1, \ldots, \Lambda_m)$ and where
$C_N^{(1)} , C_N^{(2)}$ are two positive constants which only depend
on $N$ and which will be defined in section 8. In the following, we denote $C_N^{(i)}$ some positive constant
which only depends
on $N$.

\medskip

Finally the parameter $\varepsilon >0$ is defined by
\[
\varepsilon : = p- p_N .
\]
Granted the above definitions, our result reads :
\begin{th}
\label{th6.1} Assume that $N\geq 5$ and $\mu \in {\R}$ are fixed.
Let
$(\Lambda, {\bf x} )\in P^+$ be a nondegenerate critical point of
${\cal F}_\mu $. Then, there exists $\varepsilon_0>0$ and for all
$\varepsilon \in (0, \varepsilon_0)$ there exists $u_p$ a solution
of (\ref{eq1.1}) with $\lambda := \mu \, \varepsilon^{N-4\over
N-2}$ and $p = p_N+\varepsilon$, such that
\[
|\nabla u_p|^2 \, dx  \quad \rightharpoonup \quad C_N^{(3)} \,
\sum_{i=1}^m \ell_i \, \delta_{x_i} ,
\]
in the sense of measures, where the constant $C_N^{(3)}$ is given
by
\[
C_N^{(3)} : = (N(N-2))^{{N+2}\over 4}\int_{\R^N} \left({1\over
{1+|x|^2}}\right)^{{N+2}\over 2} \, dx .
\]
\end{th}
In other words, the sequence $u_p$ converges to $0$ (in any ${\cal
C}^k$ topology) away from the points $x_i$, as the parameter $p$
tends to $p_N$. Near each $x_i$ the solution $u_p$ has multiple
blow up in the sense that there exists $c >0$ (independent of $p$), $x_{i,p}\in\O$ 
and parameters $d_{i,j,p,\mu}>0$ such that
\[
\frac{1}{c} < d_{i,j,p,\mu}  < c ,
\]
\[
x_{i,p}\to x_i,
\]
and
\[
\lim_{p\rightarrow p_N} \left\| u_p( \cdot + x_{i,p}) -
(N(N-2))^{N-2\over 4} \, \sum_{j=1}^{\ell_i} \, \left( {\bar
\varepsilon_{i,j}  \over 1 + \bar \varepsilon_{i,j} ^2 \, |  \cdot
|^2} \right)^{N-2\over 2} \right\|_{L^\infty (B_{r_0})}  =0 .
\]
Here
\[
\bar \varepsilon_{i,j} : = d_{i,j,p,\mu} \, ( \varepsilon^{{1\over
2}-j} )^{2\over N-2}
\]
and $r_0 >0$ is fixed small enough. Moreover, there is a relation
between the parameters $d_{i,1,p,\mu}$  and the parameters $\Lambda_i$ since
\[
\lim_{p\to p_N}\Lambda_i^{\frac{2}{N-2}} \, d_{i,1,p,\mu} = (N(N-2))^{1\over 2}.
\]

\medskip

We briefly describe the plan of the paper. In section 3, we give some
applications and some comments. In section 4, we recall
some well known fact about radial solutions of $\Delta u + u^p =0$
when the exponent $p$ is larger than the critical Sobolev exponent
$p_N$. In section 5 and 6 we give a new proof of existence of radial solutions.
This proof is needed just because, for
the proof of Theorem~\ref{th6.1}, we need some estimates which are
not available in \cite{DDM}.  Finally, the proof of the main
result is the content of the sections 7 to 8. This proof is based
on a gluing technic already used by  Mazzeo and Pacard \cite{MP} in a
different context.

\section{Applications and comments}

\noindent {\bf Application 1} We consider the case where $m=1$ and
$\Omega=B$ is the unit ball, we recover the result of Del Pino,
Dolbeault and Musso \cite{DDM}. Indeed, given $\ell \in\N$, we
have
$$
{\cal F}_\mu (\Lambda_1,x_1) =  {\Lambda_1^2\over
(1-|x_1|^2)^{N-2}} - \mu \, C_N^{(1)} \, \Lambda_1^{4\over N-4} +
 C_N^{(2)} \, \ell \, \log \Lambda_1
$$
It is clear that, provided the constant $\mu$ is chosen
sufficiently large, this function admits two nondegenerate
critical points which we denote by $(\Lambda_{1,1},0)$ and
$(\Lambda_{1,2},0)$. Therefore, for any $\mu$ large enough, we
find two distinct solutions of (\ref{eq1.1}).

\medskip

\noindent{\bf Application 2} Now assume that $\Omega$ is "close"
to the unit ball. Then, a standard perturbation result shows that,
for a given $\ell$ and provided $\mu$ is sufficiently large, the
function ${\cal F}_\mu$ also admits two non degenerated critical
points. This fact again guaranties the existence of two distinct
solutions of (\ref{eq1.1}).

\medskip

\noindent{\bf Application 3} We consider the case where $m=2$,
$\mu=0$, $\ell_1=\ell_2$. When  $\Omega=\R^N-\overline{B(0,1)}$
the functional ${\cal F}_0$ can be explicitely written as
$$\begin{array}{ll}
{\cal F}_0 (\Lambda_1,\Lambda_2, x_1,x_2) = &
\ds{{\Lambda_1^2\over (1-|x_1|^2)^{N-2}} + {\Lambda_2^2\over
(1-|x_2|^2)^{N-2}} + C_N^{(2)} \, \ell_1 \, \log (\Lambda_1\Lambda_2)}\\
&- \ds{2\Lambda_1\,\Lambda_2\left({1\over |x_1-x_2|^{N-2}} -
{1\over (1+ |x_1|^2|x_2|^2-2\langle
x_1,x_2\rangle)^{N-2\over2}}\right) }
\end{array}
$$
It admits a critical point $(\Lambda_1^0,\Lambda_2^0, x_1^0,x_2^0)$ where
$$
\Lambda_1^0=\Lambda_2^0=\left( C_N^{(2)} \, \ell_1\over {2\over (2a_*)^{N-2}} -{2\over (a_*^2-1)^{N-2}}-
{2\over (a_*^2+1)^{N-2}}\right)^{1/2}
$$
and
$$
 x_1^0=(a_*,0,...,0)\qquad x_2^0=(-a_*,0,...,0)
$$
where $a_* > 1$ satisfies
$$
{1\over (2a_*)^{N-1}}={a_*\over (a_*^2-1)^{N-1}}+{a_*\over (a_*^2+1)^{N-1}}.
$$
These explicit formula allow one to study (\ref{eq1.1}) for
$\mu=0$ in a annular domain $\O = B(0,1) - \overline{B(0,\rho)}$,
when $\rho$ is close to $0$. Indeed, we write $z=(z_1,z')\in
\R\times \R^{N-1}$ and using the symmetries, it is enough to look
for solutions of (\ref{eq1.1}) which only depend on $z_1$ and
$|z'|$ and blow up at two points (which turn out to be close to
$\partial B(0, \rho)$ as $p$ tends to $p_N$). For this purpose, we
study the functional ${\cal F}_0$, reduced by the symmetries we
impose. In a neighborhood of $(\rho^{N-2\over 2}\Lambda_1^0,
\rho^{N-2\over 2}\Lambda_2^0,  \rho a_*, \rho a_*)\in \R^4_+$,
${\cal F}_0$ can be expanded as
 $$\begin{array}{ll}
 {\cal F}_0 (\Lambda_1,\Lambda_2, s,t) =&  \ds{\rho^{2-N}({\Lambda_1^2\over
 (1-(s/\rho)^2)^{N-2}} +{\Lambda_2^2\over
 (1-(t/\rho)^2)^{N-2}}) }\\
  &- \ds{2\Lambda_1\,\Lambda_2\left({1\over
 (s+t)^{N-2}} -{\rho^{2-N}\over (1+ (st)/\rho^2)^{N-2}}\right) }\\
 &+  C_N^{(2)} \, \ell_1 \, \log (\Lambda_1\Lambda_2)+{\cal O}(\Lambda_1^1+\Lambda_2^2)
 \end{array}
$$
And this functional admits a non degenerated critical point,
provided $\rho$ is sufficiently small. Applying the result of
Theorem \ref{th6.1}, we find solutions of (\ref{eq1.1}) which have
two bubble trees located near $\partial B(0, \rho)$. When
$\ell_1=\ell_2=1$ such a result has been obtained by Felmer, Del
Pino and Musso in \cite{FDM} (see also \cite{FDM1} and \cite{FDM2}).

\medskip

\noindent{\bf Comments} When $m=1$, a necessary condition for the
existence of critical points of ${\cal F}_\mu $ is given by :
$\mu$ is a sufficiently large positive number. Indeed, a
nonexistence result for single peaked solutions of (\ref{eq1.1}),
when $\mu=0$, has been proved very recently by Rey {\it et al} in
\cite{Orey}.

\medskip

When $m\geq 2$ and $\mu=0$, if ${\cal F}_0 $ admits a
nondegenerate critical point, then  ${\cal F}_\mu $ also admits
  a nondegenerate critical point, provided $\mu$ is small
enough. This means that even for negative values of $\mu$, we can
construct solutions of (\ref{eq1.1}).

\medskip

Also observe that in the case where $\mu=0$, if $(\Lambda^0,{\bf
x}^0)$ is a nondegenerate critical point of ${\cal F}_0 $ for
$(\ell_1, \ldots ,\ell_m)$, then $(\sqrt{k}\Lambda^0,{\bf x}^0)$
is a nondegenerate critical point of ${\cal F}_0 $ for $(k\ell_1,
\ldots ,k\ell_m)$, where $k\in\N$.

\medskip

Finally, observe that our result parallels the corresponding
result which has been obtained by Bahri, Li and Rey \cite{BLY} for
the subcritical case, i.e. when $p < p_N$. In such case, only simple bubbles can 
be appeared, i.e. there are no bubble-towers (see \cite{Li}).

\section{Positive radial solutions of $\Delta u+u^p=0$ in $\R^N$}

We recall some well known facts about positive radial solutions of
\begin{equation}
\label{eq2.1} \Delta u+u^p=0 ,
\end{equation}
in $\R^N$. It is standard to look for radial positive solutions of
(\ref{eq2.1}) of the form
\begin{equation}
u(x) =|x|^{-\frac{2}{p-1}}v(-\log|x|). \label{eq2.2}
\end{equation}
If we set $t=-\log|x|$, then $v$ is a solution of an autonomous
second order nonlinear ordinary differential equation~:
\begin{equation}
 \partial_t^2 v  - a_p \, \partial_t v - b_p \, v + v^p = 0 ,
\label{eq2.3}
\end{equation}
where the constants $a_p$ and $b_p$ are given by \beqr
\label{eq2.4} a_p: =N-2-\frac{4}{p-1}, \label{eq2.5} \qquad
\mbox{and} \qquad b_p :=\frac{2}{p-1} \left( N-\frac{2p}{p-1}
\right). \eeqr Observe that $a_p$ vanishes precisely when $p =p_N$
and $b_p$ vanishes when $p = \frac{N}{N-2}$. We introduce the
function \beq \label{eq2.6} H_p (x,y) : = \frac{1}{2} \, y^2 -
\frac{b_p}{2} \,
 x^{2}+\frac{x^{p+1}}{p+1}.
\eeq
If $v$ is a solution of (\ref{eq2.3}), then
$$
\partial_t H_p \left( v, \partial_t v \right) =a_p \, \left(
\partial_t v \right)^2 .
$$
In particular, this implies that $\partial_t H_p \left( v,
\partial_t v \right) \geq 0$ when $p \geq p_N$.

\medskip

There are two stationary solutions of (\ref{eq2.3}), the first one
is given by $v \equiv 0$ and the other one is given by $v \equiv
b_p^{\frac{1}{p-1}}$.  We set
\[
c_p : = b_p^{\frac{1}{p-1}}.
\]
We claim that there exists a heteroclinic solutions of
(\ref{eq2.3}) when $p>p_N$. This is the content of the following~:
\begin{pro}
\label{pro2.1} Assume $p>p_N$. Then, there exists a unique
solution $v_p$ of (\ref{eq2.3}) which is defined on ${\R}$,
satisfies \beqr \label{eq2.8} \lim_{t\to-\infty} v_p(t) = c_p ,
\qquad \qquad \label{eq2.9} \lim_{t\to+\infty} v_p(t)=0 , \eeqr
and is normalized so that \beq \label{eq2.10} \lim_{t\to+\infty}
e^{\frac{2}{p-1}t}v_p(t) = 1. \eeq This solution satisfies $H(v_p,
\partial_t v_p) < 0$.
\end{pro}
{\it Proof}. We first prove that there exists a unique solution of
(\ref{eq2.3}) which is defined for $t$ large enough and which
satisfies (\ref{eq2.10}). According to a classical result in the
theory of nonlinear ordinary differential equations \cite{CL}, it
is enough to check that there exists a solution of the homogeneous
problem associated to the linearized ordinary differential
equation at $v\equiv 0$, which has the desired behavior as $t$
tends to $+\infty$. Now, the associated homogeneous problem reads
\beq \label{eq2.11}
\partial^2_t v - a_p \, \partial_t v - b_p \, v=0. \eeq
And clearly it has two
independent solutions which are given by $ t \longrightarrow
e^{\gamma_\pm t}$ where \beq \label{eq2.122} \gamma_{+} : = N-2
-\frac{2}{p-1}, \qquad \mbox{and} \qquad \gamma_{-} : =
-\frac{2}{p-1}. \eeq Therefore, there exists a unique solution of
(\ref{eq2.3}) which is asymptotic to $t \longrightarrow
e^{\gamma_- t}$ as $t$ tends to $+\infty$ and hence satisfies the
second formula of (\ref{eq2.9}). A priori this solution, which
from now on is denoted by $v_p$, is only defined for $t$ large
enough, say $t \in (\bar t,+\infty)$. Observe that there also
exists another solution of (\ref{eq2.3}) which is asymptotic to $t
\longrightarrow e^{\gamma_+ t}$ as $t$ tends to $-\infty$.

\medskip

Since the function $t \longrightarrow H_p(v_p, \partial_t v_p )$
is increasing and
$$
\lim_{t\to+\infty}{H_p \left( v_p , \partial_t v_p \right)}=0.
$$
we conclude that $H_p (v_p,\partial_t v_p )<0$ for any $t\in (\bar
t,+\infty)$. Thus, $v_p$ remains bounded independently of the
value of $\bar t$ and hence can be extended to all ${\R}$. Now, as
$t$ tend to $-\infty$, there two possibilities~: either $v_p$
converges to a limit cycle or $v_p$ converges to the constant
$c_p$, the unique stationary point in region $\{ (v,\partial_t v
 ) \, : \,  H_p (v,\partial_t v )<0\}$. But
$\partial_t H_p (v,\partial_t v ) > 0$ if $ \partial_t v \not=0$.
Hence, there are no limit cycle. We conclude that
$\lim_{t\to-\infty}v_p = c_p .$ This completes the proof of the
result. \cqfd

\medskip

In the next result, we show that the function $\partial_t v_p$
vanishes at infinitely many points, provided $p$ is close enough
to $p_N$.
\begin{pro}
\label{pro2.2} Assume that $p > p_N$ and further assume that \beq
\label{eq2.27} a_p^2-4(p-1)b_p<0. \eeq Then the set of zeros of
$\partial_t v_p $ is given by two sequences $(\underline
t_i)_{i\geq 1}$ and $(\overline t_i)_{i\geq 1}$ tending to
$-\infty$ and satisfying
\[
\overline t_i> \underline t_i > \overline{t}_{i+1}
> \underline{t}_{i+1},
\]
Moreover, we have
\[
v_p(\underline t_i) < c_p < v_p(\overline t_i).
\]
\end{pro}
{\it Proof.} We linearize (\ref{eq2.3}) at
$v=c_p$. This yields the operator
\[
L_p = \partial_t^2 - a_p \, \partial_t + (p-1) \, b_p ,
\]
since $c_p^{p-1} = b_p$. The characteristic roots of $L_p$ are
given by \beq \label{eq2.28} \tilde\gamma_\pm= \frac{1}{2} \,
(a_p\pm \, i \, \sqrt{4(p-1)b_p - a_p^2 }). \eeq These are
imaginary valued since (\ref{eq2.27}) is satisfied. It follows
from standard theory for ordinary differential equations
\cite{CL} that $v_p-c_p$ is asymptotic to a solution of the
homogeneous system associated to $L_p$. Hence there there exists
$c, d \in \R$ and $\gamma > \frac{a_p}{2}$ such that \beq
\label{eq2.29} v_p = c_p + c \, \Re \,  (e^{ \tilde \gamma_+ t
+d}) + {\cal O}(e^{\gamma t}) , \eeq as $t$ tends to $-\infty$,
where $\Re(\cdot)$ is real part of a complex number. This
immediately implies that $\partial_t v_p$ has infinitely many
zeros. The result of the proposition follows at once from this
expansion. \cqfd

\medskip

We define
\[
d_p : = \left(\frac{p+1}{2}b_p\right)^{\frac{1}{p-1}}.
\]
We now derive an upper bound for the solution
$v_p$ which has been defined in
Proposition~\ref{pro2.1}. This upper bound
follows from the more general result~:
\begin{pro}
\label{pro2.3} Assume that $v$ is a solution of
(\ref{eq2.3}) such that
$H_p (v, \partial_t v)\leq 0$ on $(t_1, t_2)$. Then $|v| \leq d_p$
on $(t_1, t_2)$.
\end{pro}
{\it Proof.} This follows at once from the fact that
$$
\max\{x > 0 \, : \, \exists \,y \in {\R} \qquad  H_p(x,y)\leq 0\}=
d_p.
$$
together with the fact that we have assumed that $H_p(v,\partial_t
v )\leq 0$.\cqfd

\medskip

From now on we assume that (\ref{eq2.27}) is satisfied and we
define the sequences
\[
\varepsilon_{p,i} = v_p (\underline t_i), \qquad \mbox{and} \qquad
\eta_{p,i} = v_p (\overline t_i) ,
\]
which correspond to the sequence of local minima and local maxima
of the function $v_p$. Observe that we have the sequence
$(\varepsilon_{p,i})_i$ (resp. $(\eta_{p,i})_i$) is increasing
(resp. decreasing) and converges to $c_p$
$$
0<\varepsilon_{p,1} <\varepsilon_{p,2} < \dots < c_p <\dots
<\eta_{p,2} < \eta_{p,1} < d_p.
$$
It will be convenient to agree that
\[
\underline t_0 =+\infty \qquad \mbox{and} \qquad \varepsilon_{p ,
0}=0.
\]
We now derive a precise expansion of the value of
$\varepsilon_{p,i}$ as $p$ tends to the critical exponent $p_N$.
This result relies on the following more general result which
gives the asymptotic of the first return map when $p$ is close to
$p_N$.

\medskip

For $p > p_N$ and $\eta \in [\varepsilon_{p,1}, c_p]$, we consider the function
$v_{p, \eta}$ which is a solution of (\ref{eq2.3}) which is
defined in $(0,\underline t_{p, \eta})$ and satisfies
\[
v_{p, \eta} (0) = \eta \qquad \mbox{and} \qquad \partial_t v_{p,
\eta}(0) = \partial_t v_{p, \eta} (\underline t_{p, \eta}) = 0.
\]
If $\underline t_{p, \eta} =+\infty$, we agree that  the above equalities
have to be understood as limits. We further assume that $v_{p,
\eta}$ is strictly increasing on $(0, \bar t_{p, \eta})$ and
strictly decreasing on $(\bar t_{p, \eta},\underline  t_{p, \eta})$. Finally,
we assume that
\[
H_p (v_{p, \eta}, \partial_t v_{p, \eta}) \leq 0 ,
\]
on $(0, \underline t_{p, \eta})$. In other words, $\underline t_{p, \eta}$ is the first
return time. Observe that, when $p=p_N$ the equation satisfied by
$v_{p, \eta}$ is Hamiltonian hence we have
\[
\lim_{p\rightarrow p_N} (v_{p, \eta}(0) - v_{p, \eta} (\underline t_{p,
\eta}))= 0.
\]
We make this estimate more precise in the following~:
\begin{pro}
\label{pro2.4} There exists  a bounded positive function $ D_N :
[0, c_{p_N}] \longrightarrow {\R}^+$, which only depends on $N$
such that
\[
\lim_{p \to p_N}  \frac{v_{p, \eta}(0)^2 -v_{p,
\eta}(\underline t_{p,\eta})^2}{p-p_N} = D_N (\eta)
\]
uniformly with respect to $\eta$.
\end{pro} {\it Proof.} For the sake of simplicity in
the notations, we drop the $p, \eta$ indices. Since  we have
assumed that $H_p (v,\partial_t v)\leq 0$ in $(0, \underline t)$, we get
$$
\left| \partial_t v\right| \leq \sqrt{b_p \, v^2-\frac{2}{p+1} \,
v^{p+1}}< \sqrt{b_p} \, v.
$$
Recall that \beq \label{eq2.33} \partial_t H_p \left(v,\partial_t
v\right)=a_p {\left(\partial_t v \right)^2}. \eeq Integrating this
equality over $(0 ,\bar t)$ and using the fact that $\partial_t v
>0$ on $(0, \bar t)$, we get
 \beq
\label{eq2.34}
\begin{array}{rllllll}
0\leq H_p (v(\bar t),0) - H_p (v(0),0) & = & \displaystyle a_p
\, \int^{\bar t}_{0}{\left(\partial_t v \right)^2 \, dt}  \\[3mm]
& \leq & \displaystyle a_p \, \sqrt{b_p} \, \int^{\bar
t}_{0}{\partial_t v \, v  \, dt} \\[3mm]
& \leq &  \displaystyle  \frac{1}{2} \,  a_p \, \sqrt{b_p} \,
v(\bar t)^2.
\end{array} \eeq
Similarly, using an integration over $(\bar t,\underline t)$ together with
the fact that $\partial_t v <0$ over this set, we also get \beq
\label{eq2.38} 0\leq H_p (v(\underline t) ,0) - H_p (v(\bar t),0)
\leq\frac{1}{2} \, a_p \, \sqrt{b_p} \,  v(\bar t)^2. \eeq Hence,
we conclude that \beq \label{eq2.39.0} 0\leq H(v(\underline t) ,0) -
H(v(0),0) \leq \, a_p \, \sqrt{b_p} \, v(\bar t)^2. \eeq Thanks to
the previous Proposition we know that $v(\bar t) \leq d_p$ and
clearly \beqr \lim_{p\to p_N}\frac{a_p}{p-p_N}= \frac{(N-2)^2}{4}
, \eeqr while $b_p$ and $d_p$ remain bounded as $p$ tends to
$p_N$. This, together with (\ref{eq2.38}) and (\ref{eq2.39.0}),
implies that \beqr \lim_{p\to p_N} \left( H_p (v(0),0) - H_p
(v(\underline t),0) \right) = \lim_{p\to p_N} \left( H_p (v(0),0) - H_p
(v(\bar t),0) \right) =0 \eeqr uniformly with respect to $\eta$.
As a consequence, we get using the expression of $H_p$ the fact
that
\[
\lim_{p \to p_N } \left( v(0)^2 - v(\underline t)^2\right) = \lim_{p \to p_N
} \left( \bar v_0 - v(\bar t)\right)  = 0
\]
uniformly with respect to $\eta$, where $\bar v_0 > v(0)$ is the
unique solution of $H_{p}(\bar v_0 , 0) = H_{p}(v(0) ,0)$ which
belongs to $(c_p, d_p)$.

\medskip

This being understood, we write
$$
\begin{array}{rlllll}
\ds \int^{\bar t}_{0}{\left(\partial_t v \right)^2dt} & = & \ds
\int_{0}^{\bar t}{\sqrt{2H_p(v (s),\partial_t v(s)) + b_p \, v^2 - {2\over p+1}
\, v^{p+1}} \, \partial_t v \, ds} \\[3mm]
 & = & \ds \int_{v(0)}^{ v(\bar t)} {\sqrt{2H_p(x,\partial_t v(x)) + b_p \, x^2 -
{2\over p+1} \, x^{p+1}} \,  dx}.
\end{array}
$$
Now, as  $p$ tends to $p_N$, it follows from the previous
discussion that the right hand side converges (uniformly with
respect to $\eta$) to
$$
E_N (\eta) : = \int_\eta^{\bar \eta} \sqrt{ 2 \, H_{p_N} (\eta, 0)
+ {(N-2)^2\over 4}x^2-{N-2\over N}x^{2N\over N-2}} dx
$$
where $ \bar \eta \geq c_{p_N}$ satisfies $H_{p_N}(\eta, 0) =
H_{p_N}(\bar \eta, 0)$. Similarly, we have
$$
\lim_{p \to p_N } \int_{\bar t}^{\underline t}{\left(\partial_t v\right)^2dt}
= E_N (\eta)
$$
where the convergence is uniform with respect to $\eta$. Moreover
the function $\eta \longrightarrow E_N (\eta)$ is bounded. Using
these limits together with (\ref{eq2.33}), which we integrate over
$(0, \underline t)$, we conclude that there exists a constant $\hat E_N
(\eta):={(N-2)^2\over 2}E_N
(\eta)$ only depending on $N$ such that \beqr \lim_{p\to p_N}
\frac{H_p(v(0),0) - H_p(v(\underline t), 0)}{p-p_N} = - \hat E_N (\eta) \eeqr
uniformly with respect to $\eta$. The result follows at once from
these limits and the expression of $H_p$.\cqfd

\medskip

Looking at the previous proof, it should be clear that
\begin{pro}
\label{pro2.5}
As $p$ tends to $p_N$, the functions
\[
\tilde v_{p, \eta} : = v_{p, \eta} (\cdot + \bar t_{p, \eta})
\]
converge (uniformly on compacts) to $w_{p_{N},\eta}$ the unique
solution of
\[
\partial_t^2 w   - b_{p_N} \, w +
w^{p_N} = 0
\]
with $w(0) = \bar \eta$ and $\partial_t w(0)=0$ where $\bar \eta \geq
c_{p_N}$ satisfy $H_{p_N} (\eta, 0) = H_{p_N} (\bar \eta, 0)$.
Moreover, the convergence is uniform with respect to $\eta$.
\end{pro}
{\it Proof.} This follows at once from Ascoli's theorem since
$v_{p, \eta}$ and all its derivatives are uniformly bounded. \cqfd

\medskip

Observe that, in the previous Proposition, as $\eta$ tends to $0$,
the function $w_{p_{N},\eta}$ converges (uniformly on compacts) to
$w_0$ which is explicitly given by
\[
w_0(t) : = \left( \frac{N(N-2)}{4} \right)^{\frac{N-2}{4}} \,
\left( \cosh t \right)^{\frac{2-N}{2}}.
\]

Going back to the study of the function $v_p$, the result of
Proposition~\ref{pro2.4} yields:
\begin{cor}
\label{pro2.8} There exists a positive constant $C_N^{(4)}$ (in
fact $C_N^{(4)}= D_N(0)$ given in Proposition~\ref{pro2.4}), only
depending on $N$, such that, for all $i \in {\N}$ \beqr \lim_{p\to
p_N}\frac{\varepsilon_{p,i}^2}{p-p_N} = i \, C_N^{(4)} \eeqr
Moreover, we have the explicit formula for $C_N^{(4)}$
\beq
C_N^{(4)}=\left( \frac{N(N-2)}{4} \right)^{\frac{N-2}{2}}{(N-2)^2\over 2(N-1)}\int_{-\infty}^{+\infty}{dt\over (\cosh t)^{N-2}}
\eeq
\end{cor}

In the next result, we estimate any solution of (\ref{eq2.3}),
near one of the points where it achieves a minimum, by comparing
it to the solution of a linear problem. Indeed, we consider $w_p$
to be the solution of the second order linear ordinary
differential equation \beq \label{eq2.40} \left\{
\begin{array}{c}
\displaystyle \partial_t^2 w_p  - a_p \, \partial_t w_p - b_p \,  w_p = 0 \\[3mm]
\displaystyle w_p(0) = 1 , \qquad \partial_t w_p(0)=0.
\end{array}
\right. \eeq  which is explicitly given by
$$
w_p = \frac{1}{N-2} \, ( \gamma_+ \, e^{\gamma_-t} - \gamma_- \,
e^{\gamma_+t}),
$$
where $\gamma_\pm$ have been defined in (\ref{eq2.122}). The
following Lemma shows that, close to $0$,  the solution $v_{p,
\eta}$ of (\ref{eq2.3}) with  $v_{p, \eta}(0)=\eta$ and
$\partial_t v_{p, \eta}(0)=0$ is well approximated by $\eta \,
w_p$.
\begin{lem}
\label{lem2.5} For all $k \in {\N}$, there exists a positive
constant $c_k >0$ such that for all $t \in {\R}$ \beq
\label{eq2.41} |\partial^k_t (v_{p, \eta} - \eta \,w_p) | \leq c_k
\, \eta^p \, w_p^p \eeq for $p$ close enough to $p_N$.
\end{lem}
{\it Proof.} Again we drop the indices $p, \eta$ to keep the
notations simple. We view $v$ as a solution of a non homogeneous
linear second order ordinary differential equation
\[
\partial_t^2 v  - a_p \, \partial_t v  - b_p \,  v
= -  v^p
\]
The variation of the constant formula yields \beq \label{eq2.42}
v(t)= \eta \,
w(t)-e^{\gamma_+t}\int_0^t{e^{(a_p-2\gamma_+)s}}\int_0^s
{e^{(-a_p+\gamma_+)\zeta} \, v(\zeta)^p \, d\zeta \, ds}, \eeq
This in particular implies that $v(t)\leq \eta \, w(t)$ for all
$t\in\R$.

\medskip

When $t \geq 0$, we can therefore use the bounds \beq v(t)\leq
\eta \, w(t) \leq c \, \eta \, e^{\gamma_+ t} \eeq in
(\ref{eq2.42}) to conclude that \beq \label{eq2.43} | v(t) - \eta
\, w(t)|\leq c \, \eta^p \, e^{\gamma_+t} \,
\int_0^t{e^{(a_p-2\gamma_+)s}} \int_0^s{e^{(-a_p+\gamma_+)\zeta}
e^{\gamma_+p\zeta} \, d\zeta \, ds} \leq c \, \eta^p \,
e^{p\gamma_+t}. \eeq When $t\leq 0$, a similar analysis yields
\beq \label{eq2.44} |v(t)- \eta \, w(t)| \leq c \, \eta^p \,
e^{p\gamma_-t} ~\mbox{for all}~t\in(-\infty,~0). \eeq This
completes the proof of the estimate of $v$. The estimates for the
derivatives follow similarly. \cqfd

\medskip

The last result translates for the function
\[
u_{p, \eta} (x)  : = |x|^{-\frac{2}{p-1}} \, v_{p, \eta} (-\log
|x|)
\]
and we obtain the estimate
\[
\left| (r \,\partial_r )^k \left( u_{p, \eta} (x)- \eta \left(
\frac{\gamma_+}{N-2} - \frac{\gamma_-}{N-2} |x|^{2-N}
\right)\right)  \right| \leq c_k
 \, \eta^p \, \left( |x|^{- p \, \gamma_-} +
|x|^{-p \, \gamma_+} \right) \, |x|^{\gamma_-},
\]
where the constant $c_k>0$ only depends on $k$ and $N$ and remains
bounded as $p\rightarrow p_N$.

\medskip

As a consequence, we have the following result which provides an
expansion of $\underline t_i$ and $\overline t_i$ as $p$ tends to
$p_N$~:
\begin{cor}
\label{cor2.5} As $p$ tends to $p_N$, we have
 \beqr \label{eq2.44.1}
\left| \overline t_i - \frac{2(i-1)}{N-2} \, \log \varepsilon
\right| + \left| \underline t_i - \frac{2i-1}{N-2} \, \log
\varepsilon \right| \leq c_i, \eeqr for some constant $c_i >0$
which only depends on $N$ and $i$. We recall that $\varepsilon =
p-p_N$.
\end{cor}
{\it Proof.} As $t$ goes from  $\underline t_{i+1}$ to $\overline
t_{i+1}$, the function $v_p$ passes once through the value $c_p$.
Hence there exists $\underline t_{*,i+1} \in (\underline t_{i+1},
\overline t_{i+1})$ such that $v ( \underline t_{*, i+1}) = c_p$.

\medskip

We first estimate $\underline t_{*,i+1} - \underline t_{i+1}$. In
view of the previous Proposition, this quantity can be estimated
by \beq \label{eq2.45} \underline t_{*,i+1}-\underline t_{i+1}=
-{1\over 2\gamma_+}\log \varepsilon + {\cal O} (1) \eeq

Now, we claim that $\overline t_{i+1} - \underline t_{*,i+1}$
remains uniformly bounded as $p$ tends to $p_N$. Indeed, it
follows from the remark after Proposition 5 that, as $p$ converges
to $p_N$, the sequence of functions $t\rightarrow v(\overline
t_{i+1}+t)$ converges on compacts to $w_0(t) = ({N(N-2)\over
4})^{\frac{N-2}{4}} \left( \cosh t \right)^{\frac{2-N}{2}}$. From
this we conclude that it takes a finite time for $w_0$ to go from
$c_{p_N}$ to $d_{p_N}$. Hence, provided $p$ remains close to
$p_N$, the time it takes to $v_p$ to go from $c_{p}$ to $v_{p}
(\overline t_{i+1})$ is bounded uniformly as $p$ tends to $p_N$.

\medskip

Therefore, we conclude that \beq \label{eq2.47} \overline
t_{i+1}-\underline t_{i+1}= - {1\over 2\gamma_+} \, \log
\varepsilon + {\cal O}(1). \eeq

\medskip

Similarly, we obtain \beq \label{eq2.48} \underline t_{i} -
\overline t_{i+1} = {1\over 2\gamma_-} \, \log \varepsilon + {\cal
O}(1). \eeq

In order to obtain the estimates as stated, just observe that
\[
{1\over \gamma_+} = \frac{2}{N-2}+ {\cal O} (\varepsilon) ,\qquad
\qquad - {1\over \gamma_-}=\frac{2}{N-2}+{\cal O} (\varepsilon ),
\]
and also that $\overline t_1={\cal O}(1)$. \cqfd

\medskip

Now, we compare solutions of (\ref{eq2.3}) which have different
boundary data. We keep the previous notations. We prove the following technical result~:

\begin{lem}
\label{lem2.6} For all $c_0>1$, there exists a positive constant $c
> 0 $ only depending on $N$ and $c_0$ such that \beqr \label{eq2.50}
{1\over c} \, (\tilde \eta -\eta) \leq |v_{p, \eta} (\underline t_{p,
\eta})- v_{p, \tilde \eta} (\underline t_{p, \tilde \eta})| \leq c \,
(\tilde \eta -\eta) \eeqr and \beqr |\underline t_{p, \eta} - \underline
t_{p, \tilde \eta}|\leq c \, \frac{\tilde \eta -\eta}{\eta} \eeqr
for all $p$ close enough to $p_N$, provided $({1\over c_0}+\sqrt{C_N^{(4)}})
\varepsilon^{1\over 2}<\eta < \tilde \eta < (c_0+\sqrt{C_N^{(4)}}) \,
\varepsilon^{1\over 2}$.
\end{lem}
{\it Proof.}  We set $v = v_{p, \eta}$ and $\tilde v = v_{p,
\tilde \eta}$. To prove the result we write for the difference $D
: = \tilde v - v$
\[
\partial_t^2 D - a_p \, \partial_t D - b_p \,
D  =  - f \, D
\]
where
\[
f : = \frac{\tilde v^p -v^p}{\tilde v -v}
\]
It follows from the estimates of Lemma \ref{lem2.5} that, for all $p$
close enough to $p_N$,
\[
|f|\leq c \, (\eta \, w_p)^{p-1}
\]
for some constant $c$ which only depends on $N$ and $c_0$. Now,  as in the proof
of Lemma \ref{lem2.5}, we use the variation of the constant formula to
get
\beq
\label{eq2.51}
D = (\tilde \eta - \eta) \, w_p - w_p \, \int_0^t e^{a_p s} \,
w_p^{-2}(s) \int_0^s e^{-a_p  \zeta} \, w_p (\zeta)
 \, f(\zeta) \, D (\zeta) \, d\zeta\, ds
\eeq
We are interested in the range of validity of the two sided
estimate
\beq
\label{eq2.52}
\frac{1}{2} \, (\tilde \eta -\eta) \, w_p \leq |D | \leq 2 \,
(\tilde \eta -\eta) \, w_p
\eeq
Inserting this into (\ref{eq2.51}), we get
\[
(\tilde \eta -\eta)  \, (1 - c \, (\eta \, w_p)^{p-1}) \, w_p \leq
|D | \leq ( \tilde \eta -\eta ) \, (1 + c \, (\eta \, w_p)^{p-1})
\, w_p
\]
Form which it follows that (\ref{eq2.52}) is valid up to the time $\hat
t_{p, \eta}$ where $c \, (\eta \, w_p)^{p-1} = 1/2$. Therefore, we
have
\beq
\label{eq2.53}
\frac{1}{2} \, (\tilde \eta -\eta)  \leq \eta \, |D (\hat t_{p,
\eta})| + \eta \, |\partial_t D (\hat t_{p, \eta})| \leq 2 \, (
\tilde \eta -\eta )
\eeq
at this point. Now, it should be clear that $
\hat t_{p, \eta}-\bar t_{p, \eta}$ is bounded independently of $\eta$ for $p$
close to $p_N$ (since $(v, \partial_t v)$ remains bounded away
from $0$ in this interval). Hence we also have
\beq
\label{eq2.54}
\frac{1}{c} \, (\tilde \eta -\eta)  \leq \eta \, |(\tilde v -v)
(\bar t_{p, \eta})|+  \eta \, |\partial_t (\tilde v  -v) (\bar
t_{p, \eta})| \leq c \, ( \tilde \eta -\eta )
\eeq
for some constant $c >0$. Standard result on dynamical systems
imply that
\beq
\label{eq2.55}
\eta \, |\bar t_{p, \tilde \eta} - \bar t_{p, \eta} | \leq c \,
(\tilde \eta - \eta).
\eeq
Using (\ref{eq2.33}) and (\ref{eq2.54}), we have
\beq
\label{eq2.56}
\begin{array}{llllll}
 &|H_p(v(\bar t_{p, \eta}),0)-H_p(\tilde v(\bar t_{p, \eta}),
\partial_t \tilde v(\bar t_{p, \eta}))| \\
 =&
 \displaystyle {|H_p(\eta,0)-H_p(\tilde \eta,0)-a_p
\int^{\bar t_{p, \eta}}_{0}{\left((\partial_t v)^2- (\partial_t \tilde v)^2\right)\, dt|}}  \\
\leq & c\eta(\tilde \eta-\eta)
\end{array}
\eeq since $({1\over c_0}+\sqrt{C_N{(4)}}) \varepsilon^{1\over 2}<\eta <\tilde \eta< (c_0+\sqrt{C_N{(4)}})
\, \varepsilon^{1\over 2}$. Together with (\ref{eq2.55}), we
estimate \beq \label{eq2.57}
 |H_p(v(\bar t_{p, \eta}),0)-H_p(\tilde v(\bar t_{p, \tilde \eta}),0)|
\leq  c\eta(\tilde \eta-\eta)
\eeq
which implies
\beq
\label{eq2.58}
 \frac{1}{c} \eta(\tilde \eta-\eta)\leq \, | \tilde v (\bar t_{p,
 \tilde \eta}) - v (\bar t_{p, \eta}) | \leq c \eta(\tilde \eta-\eta).
\eeq
From Corollary \ref{pro2.8} and results on system dynamic, there holds
\beq
\label{eq2.59}
 {1\over2 c_0} \, \varepsilon^{1\over 2}<v(\underline  t_{p, \eta})<
\tilde v(\underline  t_{p,\tilde \eta}) <2c_0 \, \varepsilon^{1\over 2}.
\eeq Using similar arguments on $(\bar t_{p, \eta}, \underline t_{p,
\eta})$, we get
\[
\eta \, | (\underline t_{p, \tilde \eta} - \bar t_{p, \tilde
\eta}) - ( \underline t_{p, \eta} - \bar t_{p, \eta}) | \leq c \,
(\tilde \eta - \eta)
\]
and also that
\[
\frac{1}{c} (\tilde \eta-\eta)\leq \, | \tilde v (\underline t_{p,
 \tilde \eta}) - v (\underline t_{p, \eta}) | \leq c (\tilde \eta-\eta).
\]
The result follows at once from these estimates. \cqfd

\section{Linear results}

We keep the notations in the previous section. For the sake of simplicity in the notations, we drop the
indices $p$ and $\eta$. We consider $w$ to be the solution
of
\beq
\label{eq2.73.0}
\partial_t^2 w - a_p \, \partial_t w - b_p \, w  + p \, v^{p-1} \, w  =
 e^{-2t} \, v
\eeq in  $(0,\underline t)$ with boundary conditions $w
(\underline t)= \partial_t w(\underline t) = 0$.  We are
interested in the behavior of $w$ as $\eta$ tends to $0$. This is
the contain of the following result.
\begin{lem}
\label{lem3.1} Assume that $N \geq 5$. Let $c_0 > 1$ and $d_0 >0$
be fixed. Assume that $\eta \in ({1\over c_0} \varepsilon^{1\over
2},c_0 \varepsilon^{1\over 2})$. Then, there exist  $c >0$ and
 $\varepsilon_0 >0$ such that for all $p\in (p_N,
 p_N+\varepsilon_0)$ we have
\beq \label{eq2.73} \left| w(t)- {4 \beta_{p, \eta} \over \eta  \,
(N-2)^2} \, e^{-(N-2)t/2}\right| \leq c \, \beta_{p, \eta} \,
\varepsilon^{-{1\over 2}+{2\over N+2}} \eeq in $(-d_0, d_0)$,
where the constant $\beta_{p, \eta}$ is given by \beq
\label{eq2.73.1} \beta_{p, \eta} : = \int^{\underline t_{p, \eta}
}_0 v^2_{p, \eta} (s) \, e^{-2s} \, ds. \eeq Moreover, we have
\beq \label{eq2.7333} \lim_{p\to p_N} e^{2\bar t_{p, \eta}} \,
\beta_{p, \eta} = \left( {N(N-2)\over 4} \right)^{N-2\over 2}
\int_{-\infty}^{+\infty} {e^{-2s}\over (\cosh s)^{N-2}} \, ds\,:=C_N^{(5)}.
\eeq
\end{lem}
{\it Proof.} As usual we drop the $p, \eta$ indices. We use the
fact that $w_1 : = \eta^{-1} \,
\partial_tv$ is an explicit solution of the homogeneous problem
\beq
\label{eq2.73.2}
\partial_t^2 w_1 - a_p \, \partial_t w_1 - b_p \, w_1  + p \, v^{p-1} \,
w_1 =0.
\eeq
This yields a representation formula for $w$, at least when $t \in
(\overline t , \underline t)$.
\[
\label{eq2.74}
w (t) = w_1 (t) \int_t^{\underline t} e^{a_p s} \, w_1^{-2}(s) \,
\int_s^{\underline t} e^{- a_p \zeta} \, w_1(\zeta) \, e^{-2\zeta}
\, v(\zeta) \, d\zeta \, ds.
\]
Observe that the result of Lemma \ref{lem2.5} yields
\[
\label{eq2.75}
\frac{1}{c} \, e^{\gamma_- (t - \overline t)}  \leq v(t) \leq c
\, e^{\gamma_- (t - \overline t)}
\]
for all $t \in (\overline t, \underline t)$ and also
\[
\label{eq2.76}
\frac{1}{c} \, e^{\gamma_- (t - \overline t)}  \leq |\partial_t
v(t)| \leq c \, e^{\gamma_- (t - \overline t)}
\]
for all $t \in (\overline t + 1 , \underline t -1)$. Using this,
we get the estimate
\beq
\label{eq2.77}
|w(t)|+ |\partial_t w(t)|\leq c \, e^{-2t}\, e^{\gamma_-
(t-\overline t)}
\eeq
for some constant $c >0$. This estimate is valid for all $t\in
(\overline t +1, \underline t )$, however, enlarging the value of
$c$ if this is necessary, we can assume that this estimate holds
for $t \in (\overline t - 1, \underline t )$. The solution $w$
extends to $(0, \underline t)$.

\medskip

Again, we use the fact that $w_1 = \eta^{-1} \, \partial_t v$ and
\[
\label{eq2.78} w_2 (t) :  =  w_1 (t) \int_t^{\overline t -1}
e^{a_p s} \, w_1^{-2}(s) \, ds
\]
which is defined for $t \in (1, \overline t -1)$, are solutions of
the homogeneous problem (\ref{eq2.73.2}). Hence we can decompose
\[
\label{eq2.79} w = \alpha_1 \, w_1  + \alpha_2 \, w_2  + \tilde w
\]
where $\tilde w$ is defined by
\[
\label{eq2.80} \tilde w (t) : =w_1 (t) \int_0^{t} e^{a_p s} \,
w_1^{-2}(s) \, \int_0^{s} e^{- a_p \zeta} \, w_1(\zeta) \,
e^{-2\zeta} \, v(\zeta) \, d\zeta \, ds.
\]
As above, the result of Lemma \ref{lem2.5} yields
\beq
\label{eq2.81} \frac{1}{c} \, \eta \, e^{\gamma_+ t}  \leq v(t)
\leq c \, \eta \, e^{\gamma_+ t}
\eeq
for all $t \in (0, \overline t)$ and also \beq \label{eq2.82}
\frac{1}{c} \, \eta \, e^{\gamma_+ t}  \leq |\partial_t v(t)| \leq
c \, \eta \, e^{\gamma_+ t} \eeq for all $t \in (1, \overline
t-1)$. Using this, we get the estimate \beq \label{eq2.83} |\tilde
w(t)|+ |\partial_t \tilde w(t)|\leq c \, \eta \, e^{-2t}\,
e^{\gamma_+ t} \eeq for some constant $c >0$. This estimate is
valid for all $t\in (0, \overline t -1)$.
\medskip

Since
\[
\frac{1}{c} \, \leq \eta \, | w_1 ( \overline t -1)|
\]
it follows at once from (\ref{eq2.77}) and (\ref{eq2.83}) that we can estimate the
parameter $\alpha_1$ by
\[
\label{eq2.84}
|\alpha_1| \leq c \, \eta e^{-2\bar t}.
\]

In order to estimate the parameter $\alpha_2$ we multiply the
equation (\ref{eq2.73.0}) by $w_1$ and integrate by parts. Using
the fact that $w_1$ is a solution of (\ref{eq2.73.2}) we obtain
\[
\label{eq2.85} \left[ w_1 \, \partial_t w - w\, \partial_t w_1  -
a_p \, w\, w_1 \right]_0^{\underline t} = -2a_p  \,
\int_0^{\underline t}  \partial_t w_1 \, w dt + \eta^{-1} \,
\int_0^{\underline t} e^{-2t} \, v \, \partial_t v \, dt.
\]
Since $w_1=0$ at $t=0$ and $w = \partial_t w =0$ at $t=\underline
t$, this simplifies into
\[
\label{eq2.86} w (0) \, \partial_t w_1 (0)   = \ -2a_p  \,
\int_0^{\underline t}  \partial_t w_1 \, w \, dt +  \eta^{-1} \,
\int_0^{\underline t} e^{-2t} \, v \,\partial_t v \, dt
\]
From (\ref{eq2.81}) and (\ref{eq2.82}), it follows
\[
\label{eq2.87} |w_2(t)|+ |\partial_t  w_2(t)|\leq c \,
 e^{-\gamma_+ t}
\]
for all $t \in (1, \overline t-1)$. Enlarging the value of
$c$ if this is necessary, we can assume that this estimate holds
for $t \in (0, \overline t )$.

\medskip

Collecting these estimates, we get
\[
\label{eq2.88} \int_0^{\overline t}  w \, \partial_t w_1  \, dt =
\alpha_2  \, {\cal O}(\log \varepsilon ) + {\cal O}(
{\varepsilon}^{-{1\over 2}+{2\over N-2}}), \qquad \mbox{and}
\qquad \int_{\overline t}^{\underline t}  w\partial_t w_1 ={\cal
O}( {\varepsilon}^{-{1\over 2}+{2\over N-2}}).
\]
To calculate $w_2(0)$, we see that
\[
\label{eq2.90} {e^{a_pt}(pv^{p-1}(t)-b_p)\over (\partial_t^2
v(t))^2} ={1\over \partial_t v(t)}{d\over dt}\left({e^{a_pt}\over
\partial_t^2 v(t)}\right).
\]
Hence we get
\[
\label{eq2.90.1}
\begin{array}{lll}
w_2(t)&=& \ds{w_1(t)\int_1^{\bar t-1}e^{a_p s}w_1^{-2}(s) ds + w_1(t) \int_t^1 e^{a_p s} w_1^{-2}(s) ds}\\
&=&\ds{w_1(t)\int_1^{\bar t-1}e^{a_p s}w_1^{-2}(s) ds +w_1(t)\left[-{e^{a_p s}\over w_1(s)\partial_t
w_1(s)}\right]_t^1}\\
&&\ds{- w_1(t)\int_t^1 {e^{a_p s}(b_p-pv^{p-1}(s))\over (\partial_t w_1(s))^2}ds}
\end{array}
\]
for all $t\in (-d_0,d_0)$. In particular
\[
\label{eq2.91} w_2(0)={1\over b_p-\eta^{p-1}},
\]
since $w_1(0)=0$. Consequently, we obtain the estimate
$$
\alpha_2= \eta^{-1} \int^{\underline t}_0 v^2(s) \, e^{-2s} \, ds
+{\cal O}( {\varepsilon}^{1\over 2}).
$$

\noindent It remains to estimate $w_2$ in the neighborhood of $0$.
We first estimate $\ds{\int^{\bar t-1}_1 {e^{a_ps}ds\over
(w_1(s))^2}}$. We decompose
\[
(1,\bar t-1)= (1,{4\over N^2-4}\log{1\over \varepsilon})\cup
({4\over N^2-4}\log{1\over \varepsilon}, \bar t-1):=I_1\cup I_2.
\]
It follows from (\ref{eq2.82}) that
\beq
\label{eq2.92} \int_{I_2}{e^{a_ps}ds\over (w_1(s))^2}\leq c \,
\varepsilon^{4\over N+2} \qquad \mbox{and} \qquad \label{eq2.93}
1\leq e^{a_ps}\leq 1 + c \, \varepsilon \, \log{1\over
\varepsilon},
\eeq
for all $s \in (0, \bar t)$. Using the result of
Lemma~\ref{lem2.5}, we obtain \[ \label{eq2.95} \partial_t
v(t)={N-2\over 2} \eta \, \mbox{sinh}({N-2\over 2}t)(1+ {\cal
O}(\varepsilon^{2\over N+2})) \] for all $ t\in I_1$. Therefore,
we deduce
$$
\int_{I_1} {e^{a_ps}ds\over (w_1(s))^2}={8\over
(N-2)^3}\left({\mbox{cosh}({N-2\over 2})\over \mbox{sinh}
({N-2\over 2})}-1\right)+ {\cal O}(\varepsilon^{2\over N+2}).
$$
On the other hand, using again Lemma \ref{lem2.5}, we have
\beqr
\label{eq2.94} w_1(t)= {N-2\over 2}  \, \mbox{sinh}({N-2\over
2}t)+ {\cal O}(\varepsilon^{p-1\over
2})\\[3mm]
\label{eq2.94.1}
\partial_t w_1(t)={(N-2)^2\over 4} \,\mbox{cosh}({N-2\over 2}t) + {\cal O}(\varepsilon^{p-1\over 2})
\eeqr for all $t\in (-d_0,d_0)$. Now direct calculations lead to
$$
w_2(t)= {4\over (N-2)^2}\, e^{-{(N-2)t\over 2}}+ {\cal
O}(\varepsilon^{2\over N+2})
$$
for all $t\in (-d_0,-d_0)$. This proves (\ref{eq2.73}).

\medskip

Finally, in order to obtain (\ref{eq2.7333}), it is enough to
observe that $v(\bar t+\cdot)$ converges (uniformly on compacts)
to $w_0$. This completes the proof of the result.\cqfd

\medskip

Using similar arguments (and the notations of the previous
Proposition), one can show that \beq \label{eq2.96} \left| w(t)-
{4 \, \beta \over\eta (N-2)^2} \, e^{-(N-2)t/2} \right| \leq c \,
\beta \, \varepsilon^{-{1\over 2}+{3\over N+2}} \eeq if $N\geq 6$,
and \beq \label{eq2.97} \left| w(t) - {4 \, \beta  \over \eta
(N-2)^2} \, e^{-(N-2)t/2} \right|\leq c \beta\varepsilon^{-{1\over
2}+{3\over N+2}-{1\over N^2-4}} \eeq if $N=5$, for all $t\in
({2\over N^2-4}\log{1\over \varepsilon} - d_0,{2\over
N^2-4}\log{1\over \varepsilon}+d_0)$.

\medskip

In the following, $\beta_{p,\eta}$ will be expanded.

\begin{lem}
\label{lem3.2} Under the above assumptions, let $c_0>1$ be given. Assume
$$
{1\over c_0}\, \varepsilon^{1\over 2}\leq \eta\leq c_0\,
\varepsilon^{1\over 2}.
$$
Then, \beq \label{eq2.103} \bar t= {2\over N-2}\log{1\over
\eta}+C_N^{(6)}+O(\varepsilon^{2\over N} \log{1\over \varepsilon})
\eeq where $$C_N^{(6)}={2\over N-2}\log 2+{1\over 2}\log N(N-2).$$
\end{lem}
{\it Proof.} We recall
\beq
\label{eq2.104}
 \bar t =\int_{\eta}^{v(\bar t)}{dv\over\sqrt{2H_p(v,\partial_t v)+b_pv^2- {2v^{p+1}\over
p+1}} }
\eeq We divide
$$
\begin{array}{lll}
[\eta, v(\bar t)]&=&[\eta, \varepsilon^{N-2\over 2N}] \cup
[\varepsilon^{N-2\over 2N}, {1\over 2}({N(N-2)\over 4})^{N-2\over
4}]\cup [{1\over 2}({N(N-2)\over 4})^{N-2\over 4},
v(\bar t)]\\
&:=&I_1\cup I_2\cup I_3
\end{array}
$$
We estimate
$$
\begin{array}{lll}
\ds{\int_{I_1}{dv\over\sqrt{2H_p(v,\partial_t v)+b_pv^2-
{2v^{p+1}\over p+1}} }}&=\ds{(1+{\cal O}(\varepsilon^{2\over N}))
\int_{I_1}{dv\over
\sqrt{b_p(v^2-{\eta}^2)}}}\\
&=\ds{(1+{\cal O}(\varepsilon^{2\over N})) {2\over
N-2}Argch(\varepsilon^{N-2\over 2N}/\eta)}
\end{array}
$$
$$
\begin{array}{lll}
\ds{\int_{I_2}{dv\over\sqrt{2H_p(v,\partial_t v)+b_pv^2-
{2v^{p+1}\over p+1}} }}&=\ds{(1+{\cal O}(\varepsilon^{2\over N}))
\int_{I_2}{dv\over \sqrt{({N-2\over 2})^2v^2-{N-2\over
N}v^{2N\over N-2}}}}
\end{array}
$$
$$
\begin{array}{lll}
\ds{\int_{I_3}{dv\over\sqrt{2H_p(v,\partial_t v)+b_pv^2-
{2v^{p+1}\over p+1}} }}&=\ds{\int_{I_3}{dv\over \sqrt{({N-2\over
2})^2v^2-{N-2\over N}v^{2N\over N-2}}}+ {\cal
O}(\varepsilon^{1\over 2})}
\end{array}
$$
Recall $w_0(t)
=({N(N-2)\over 4})^{\frac{N-2}{4}} \left( \cosh t
\right)^{\frac{2-N}{2}}$. We deduce
$$
\int_{I_2\cup I_3}{dv\over \sqrt{({N-2\over
2})^2v^2-{N-2\over N}v^{2N\over N-2}}}=\tilde t
$$
where $w_0(\tilde t)=\varepsilon^{N-2\over 2N}$.
Hence, the desired results yield.
\cqfd

We set $\bar v(\cdot)=v(\bar t+\cdot)$. We have the following result

\begin{lem}
\label{lem3.3} Given $c_0>1$, assume $\bar v(0)\in
(d_p-c_0\varepsilon,d_p)$. Then, there exists  the constant $c$
independent of  $p$ such that \beq \label{eq2.105} |\bar
v(t)-w_0(t)|+ |\partial_t \bar v(t)-\partial_t w_0(t)|\leq c\,
\varepsilon\,e^{(N-3/2)|t|\over 2} \eeq for all $t\in (-\bar
t,\underline t-\bar t)$, provided $p$ close to $p_N$.
\end{lem}
{\it Proof.} We write for the difference $\tilde D:=\bar v-w_0$ so that
\beq
\label{eq2.106}
\partial_t^2 \tilde D -b_{p_N } \tilde D =-f  \tilde D +g
\eeq
where
$$
f:={\bar v^p-w_0^p\over \bar v-w_0},~~~\mbox{and}~~~g=a_p\partial_t \bar v+w_0^{p_N}-w_0^p.
$$
Clearly, there exist some positive constants $K$ and $c$ independent of $p$ such  that
$$
| -b_{p_N }+f(t)|\leq ({N-3/2\over 2})^2
$$
for all $t\in (K,\underline t-\bar t)\cup ( -\bar t,-K)$ and
$$
| g(t)|\leq  c\, \varepsilon
$$
for all $t\in ( -\bar t,\underline t-\bar t)$. Recall
$$
|\bar v(0)-w_0(0)|+ |\partial_t \bar v(0)-\partial_t w_0(0)|\leq
c_0\, \varepsilon.
$$
Hence, the desired result follows from the standard ordinary
differential equation theory. \cqfd

As a consequence, we obtain immediately

\begin{cor}
\label{cor3.1} There exists a positive constant $C_N^{(7)}$ (only
depending on $N$), such that
$$
H_p(\varepsilon_{p,\ell },0)= - \ell \, C_N^{(7)} \, \varepsilon+
{\cal O}(\varepsilon^{2-{1\over 2(N-2)}})$$
where $ C_N^{(7)}:={(N-2)^2 C_N^{(4)}\over 8}$. In particular, \beq
\label{eq2.199} \varepsilon_{p,\ell} =  \, (\ell
\varepsilon C_N^{(4)})^{1/2} + {\cal O}(\varepsilon^{N+2\over
 2N}).
\eeq
\end{cor}

We keep the notations introduced in section 4 and we define for
all $\ell \in\N$ \beq \beta_{p,\ell} := \int^{\underline
t_{\ell-1}}_{\underline t_\ell} v^2_p(t) e^{-2(t-\underline
t_\ell)} \, dt. \eeq Thanks to Lemma \ref{lem3.1} to \ref{lem3.3}
and Corollary \ref{cor3.1}, we conclude
\begin{cor}
\label{cor3.2} There exists a positive constant $C_N^{(8)}$ (only
depending on $N$) such that \beq \label{eq2.200}
\beta_{p,\ell}=(\ell \varepsilon)^{2\over N-2} \left( C_N^{(8)}+
{\cal O}( \varepsilon^{2\over N}\log{1\over
\varepsilon}+\varepsilon^{N-17/4\over N-2}) \right) \eeq
where $C_N^{(8)}:=(C_N^{(4)})^{2\over N-2}C_N^{(5)}e^{-2C_N^{(6)}}$.
\end{cor}

\section{Radial solutions of $\Delta
u + \lambda u + |u|^{p-1}u=0$ in the unit ball}

In this section we recover part of the result of Del Pino,
Dolbeault and Musso concerning the existence of solutions of (\ref{eq1.1})
in the unit ball. In doing so our aim is to derive precise
estimates for these solutions which will be needed in the
forthcoming construction.

\medskip

We begin with the definition of weighted spaces in cylindrical
coordinates. These spaces are at the heart of our construction.
\begin{defi}
Given $\delta \in\R$ and $-\infty \leq t_1 < t_2 \leq +\infty $,
the space ${\cal C}^0_\delta((t_1, t_2)\times S^{N-1})$ is defined
to be the set of continuous functions $w \in {\cal
C}^0_{loc}((t_1,t_2)\times S^{N-1})$ for which the following norm
is finite~: \beq \label{eq3.1} \|w\|_{{\cal C}^0_\delta
((t_1,t_2)\times S^{N-1})}:= \| e^{- \delta s} \, w \|_{L^\infty
((t_1,t_2)\times S^{N-1})} . \eeq
\end{defi}

We would like to prove the existence of radial solutions of
\beq
\label{eq5.2}
\begin{array}{ll}
\Delta u +\lambda u + |u|^{p-1}u =0 &\mbox{in } B(0,1)
\end{array}
\eeq Using (\ref{eq2.2}), we
reduce the study of (\ref{eq5.2}) to study of the nonlinear second
order ordinary differential equation
\beq \label{eq5.3}
\partial_t^2 v - a_p \, \partial_t v - b_p \, v  + |v|^{p-1}v +\lambda \,
e^{-2t} \, v = 0, \eeq in $(0,+\infty)$. We keep the
notations introduced in section 4 and we consider the linear
operator \beq \label{eq5.1.1} L_{p, \eta} : = \partial_t^2 - a_p
\, \partial_t  - b_p + p \, v_{p,\eta}^{p-1}. \eeq We state,
without a proof a result which will be proven in a more general
context in the next section.
\begin{pro}
\label{lem5.1} Assume that $ \delta \in (-{N+2\over 2} ,-
{N-2\over 2})$ is fixed. Then, there exist $\varepsilon_0 >0$,
$\eta_0 >0$ and $c >0$ such that, for all $\varepsilon \in (0 ,
\varepsilon_0)$, for all $\eta \in (0, \eta_0)$ and for all $f \in
{\cal C}^0_\delta ((0,\underline t_{p,\eta}))$, there exists a
unique solution $w \in {\cal C}^0_\delta ((0,\underline
t_{p,\eta}))$ of \beq \label{eq5.1.2} L_{p, \eta} \, w =  f , \eeq
in $(0,\underline t_{p,\eta})$ which satisfies \beq w(\underline
t_{p,\eta}) =
\partial_t w(\underline t_{p,\eta})=0,
\label{zeze} \eeq with $p = p_N + \e$. Furthermore, \beq
\label{eq5.1.3} \|w\|_{{\cal C}^0_\delta} \leq c \,
\|f\|_{C^0_\delta} .\eeq
\end{pro}
When $\underline t_{p, \eta} < +\infty$, the existence and
uniqueness of the solution of (\ref{eq5.1.2}) is straightforward
but the uniform estimate (\ref{eq5.1.3}) requires some work. When
$\underline t_{p, \eta} =+\infty$, the boundary data (\ref{zeze})
have to be understood as limits as $\underline t_{p,\eta}
=+\infty$.

\medskip

The next result will allow us to recover  (part of) the result of
Del Pino, Dolbeault and Musso \cite{DDM}~:
\begin{pro}
\label{th4.1} Assume that $\ell \in\N$ is fixed and that $N\geq
5$.  Then, there exists $ \varepsilon_0 > 0$ such that for all
$\mu \in\R$, for all $\xi \in \R$ and for all $\varepsilon \in (0,
\varepsilon_0)$, problem (\ref{eq5.2}) with $p= p_N+\varepsilon$
and $\lambda = \mu \, \varepsilon^{N-4\over N-2}$ admits a
solution which can be written in the form \beq \label{eq5.3.0}
u_{p, \lambda, \xi} (x) = ( N \, (N-2))^{N-2\over 4} \,
\sum_{j=1}^\ell \left( {\bar \varepsilon_j \over 1 + \bar
\varepsilon_j^2 \, |x|^2} \right)^{N-2\over 2} + o(1) \eeq where
$o(1)$ converges uniformly to $0$ on $B(0,1)$ as $\e$ tends to $0$
and where
\[
\bar \e_j : = d_j \, (\varepsilon^{{1\over 2}-j})^{2\over N-2}
\]
for some  parameters $d_j$ which are bounded from below and from
above by some positive constant independent of $\e$. Moreover we
have the following expansion
\beq
\label{eq5.1.4}
\begin{array}{ll}
&u_{p, \lambda, \xi} (x)\\
=&\ds{(\ell\varepsilon)^{1\over 2}\left[{\sqrt{ C_N^{(4)}}\over 2}e^{(N-2)\xi\over 2}
+{\sqrt{ C_N^{(4)}}\over 2}e^{(2-N)\xi\over 2}|x|^{2-N}-{4\mu  C_N^{(8)}\ell^{4-N\over N-2}\over
(N-2)^2\sqrt{ C_N^{(4)}}}e^{(N-6)\xi\over 2}
\right]}\\
&+{\cal O}(\varepsilon^{1\over 2}r_\varepsilon^2)
\end{array}
\eeq
in $B(0,2r_\varepsilon)- B(0, r_\varepsilon/2)$, where $r_\varepsilon:=\varepsilon^{2\over N^2-4}$.
Furthermore, $u_{p, \lambda, \xi}$ is positive in $B(0,2r_\varepsilon)$.
\end{pro}
{\it Proof.} The proof is decomposed in several steps. We give the
prove in the case where $N\geq 6$ since, when $N=5$, the proof is
similar with straightforward changes. Given $\xi \in {\R}$ (which
will be fixed later on) we define
$$
T_{2i}=\underline{t}_{\ell-i}-\underline{t}_\ell + \xi , \qquad
\qquad T_{2i-1} = \overline{t}_{\ell-(i-1)}-\underline{t}_\ell +
\xi
$$
for$0 < i < \ell$ and
$$
T_{2\ell}=+\infty, \qquad
T_{2\ell-1}=\overline{t}_1-\underline{t}_\ell + \xi, \qquad T_0=0,
$$
For all $0\leq i\leq \ell - 2$, we define $v_{p,i}$ to be the
solution of (\ref{eq2.3}) in $[T_{2i}, T_{2i+2}]$ with boundary
conditions
$$
v_{p,i}(T_{2i+2}) = \varepsilon_{p,\ell - 1 - i } + \alpha_i
,\qquad \qquad \partial_t v_{p,i}(T_{2i+2})=0
$$
and we define $v_{p,\ell - 1}$ to be the solution of (\ref{eq2.3})
in $[T_{2\ell-2}, +\infty)$ with boundary conditions
$$
v_{p,\ell-1}(T_{2\ell-2}) = \varepsilon_{p,1}, \qquad \qquad
\partial_t v_{p,\ell-1}(T_{2\ell-2})=0 ,
$$
for some parameters $\alpha_i \in {\R} $ (which are assumed to be
small). Here the parameters $\varepsilon_{p,i}$ are the one which
have been introduced in section 4.

\medskip

For any $0\leq i\leq \ell - 1$, we define the function \beq
\label{eq5.3.1} W_i(t) :=v_{p,i}(t+t_i)+w_i(t) \eeq for on the
interval $[T_{2i},T_{2i+2}]$, for some parameters $t_i\in \R$ and
some functions $w_i\in C^0([T_{2i}, T_{2i+2}])$. We agree that
$t_{\ell-1}=0$ and $\alpha_{\ell-1}=0$.

\medskip

Granted the above definitions, our strategy is the following~: In
Step 1 and 2, we look $W_i$ solutions of (\ref{eq5.3}) on each
interval $[T_{2i}, T_{2i+2}]$. Moreover, $W_i$ are positive if $i\geq 1$. In Step 3, we choose the parameters
$(\alpha_0,\ldots ,\alpha_{\ell-2})$ and $(t_0, \ldots
,t_{\ell-2})$ so that the Cauchy data of $W_i$ and of $W_{i-1}$
coincide at $T_{2i}$. Gathering the functions $W_i$ together, we
obtain a solution of (\ref{eq5.3}) which still depends on $\xi$.

\medskip

\noindent {\bf Step 1.} For each $1\leq i\leq \ell - 1$, we look
for a solution of (\ref{eq5.3}) in $[T_{2i},T_{2i+2}]$. Recall
that $\varepsilon=p-p_N$. We now assume that \beq \label{eq5.001}
\alpha_i=o(\varepsilon^{1/2}) \qquad  \xi =  {\cal O}(1) \qquad \mbox{and}
\qquad t_i=o(1) \eeq as $\varepsilon$ tends to $0$. We
define the operator \beq \label{eq5.4} L_{p,i} = \partial_t^2 -
a_p  \,
\partial_t  - b_p + p \, v_{p,i}^{p-1} (\cdot + t_i).
\eeq With these notations, the equation we need to solve reads
\beq \label{prpr}
L_{p,i} w_i = -\lambda \, e^{-2t}( v_{p,i} (
\cdot + t_i) + w_i ) - Q_i(w_i) \eeq where we have defined
$$
Q_i(w_i) : = |v_{p,i}(\cdot +  t_i)+ w_i|^{p-1}(v_{p,i}(\cdot +  t_i)+ w_i)- v_{p,i}^p( \cdot
+t_i) - p \, v^{p-1}_{p,i}( \cdot + t_i) \, w_i.
$$

We fix the weight parameter $\delta \in (-{N-1\over 2},-{N-2\over
2})$ and we consider the set of functions
$$
{\cal E}_{\kappa ,i} =  \left\{ w\in {\cal C}^0_\delta((T_{2i},
T_{2i+2})) \quad  : \quad  \|w\|_{{\cal C}^0_\delta } \leq \kappa
\, \lambda \, e^{-(\delta+2)T_{2i+1}} \right\},
$$
where the constant $\kappa >0$  will be fixed later on.

\medskip

Given $w\in {\cal E}_{\kappa,i}$, it follows from (\ref{eq2.44.1})
that $|w |\leq c \,  \kappa \, \lambda \,
\varepsilon^{\delta+6\over N-2}$. Recall that $\lambda = \mu \,
\varepsilon^{N-4\over N-2}$ and $\delta > - {{N-1}\over 2}$,
hence, we obtain
\[
\label{eq5.5} |w |\leq  c \, \kappa \, \varepsilon^{2+\delta
+N\over N-2} \ll v_{p,i}( \cdot + t_i)
\]
on $(T_{2i},T_{2i+2})$. Therefore, we are allowed to use Taylor's
expansion $|(1+t)^p - 1 - p\, t |\leq c \, t^2$ for $t$ close
enough to $0$, to estimate
\[
\label{eq5.6}
\begin{array}{ll}
|Q_i(w)| \leq c \, v_{p,i}^{p-2} ( \cdot + t_i) \, w^2
\end{array}
\]
Using this, we obtain
\beq \label{eq5.7}
\begin{array}{ll}
\|Q_i(w)\|_{{\cal C}^0_\delta} \leq c \,  \kappa^p  \, \lambda \,
\varepsilon^{(\delta+2+N)(p-1)\over N-2}e^{-(\delta+2) T_{2i+1}}.
\end{array}
\eeq Next, we estimate \beq \label{eq5.8} \|\lambda
e^{-2t}v_{p,i}( \cdot + t_i)\|_{{\cal C}^0_\delta} =
\ds{\sup_{(T_{2i},T_{2i+2})} \lambda \, e^{-(\delta+2)t} \,
|v_{p,i}(\cdot + t_i) |} \leq c \,  \lambda \,
e^{-(\delta+2)T_{2i+1}}, \eeq since $\gamma_+-\delta-2>0$ and
$\gamma_- - \delta - 2 < 0$ provided $\varepsilon$ is close enough
to $0$. With similar arguments, we get \beq \label{eq5.9}
\|\lambda e^{-2t} \, w  \|_{{\cal C}^0_\delta} \leq \lambda
e^{-2T_{2i}}\| w \|_{{\cal C}^0_\delta} \leq  c \, \kappa \,
\lambda^2 \, \varepsilon^{4\over N-2}\, e^{-(\delta+2) T_{2i+1}} ,
\eeq Combining (\ref{eq5.7}) to (\ref{eq5.9}), we have obtained
\beq \label{eq5.10} \| Q_i(w) + \lambda \, e^{-2t} \, (v_{p,i}(
\cdot + t_i)+ w) \|_{{\cal C}^0_\delta} \leq  c \, \lambda \,
e^{-(\delta+2)T_{2i+1}} \, ( \kappa^p \,
 \varepsilon^{(2+\delta + N)(p-1)\over N-2} + 1 + \kappa \, \lambda
\, \varepsilon^{4\over N-2}), \eeq which holds for all $w\in {\cal
E}_{\kappa,i}$.

\medskip

Given $w\in {\cal E}_{\kappa,i}$, we apply the result of
Proposition~\ref{lem5.1} which provides a solution of
\[
\label{eq5.12} L_{p,i} \, \tilde w  = - Q_i(w) -
\lambda \, e^{-2t}  \, (v_p( \cdot + t_i)+ w)
\]
with $ \tilde w(T_{2i+2}) = \partial_t \tilde w (T_{2i+2})=0$.
Thanks to (\ref{eq5.10}), we also have the estimate
\[
\label{eq5.11} \|\tilde w \|_{{\cal C}^0_\delta} \leq \lambda \,
e^{(-\delta-2)T_{2i+1}} \,  \tilde c \, (\kappa^p \,
\varepsilon^{(2+\delta +N)(p-1)\over N-2} + 1 + \kappa \, \lambda
\, \varepsilon^{4\over N-2}).
\]
for some constant $\tilde c >0$ which does not depend on $w$, nor
on $\kappa$ nor on $\varepsilon$ provided this later is chosen
small enough. This estimate being understood, we choose the
constant $\kappa
> 0$ so that
$$
\tilde c \, (\kappa^p \, \varepsilon^{(2+\delta +N)(p-1)\over N-2}
+ 1 + \kappa \, \lambda \, \varepsilon^{4\over N-2}) \leq \kappa ,
$$
for all $\varepsilon$ close enough to $0$, say $\varepsilon \in
(0, \varepsilon_0)$.

\medskip

To summarize, using the above analysis, we can define the mapping
$$
{\cal T}_{i} :  {\cal E}_{\kappa,i} \longrightarrow  {\cal
E}_{\kappa,i}
$$
by ${\cal T}_{i}(w) = \tilde w$. Thanks to the above choice of
$\kappa$, the mapping ${\cal T}_{i}$ is well defined. Observe that
this mapping is clearly continuous and compact so that one can
refer to Schauder's fixed point Theorem to obtain the fixed point
of ${\cal T}_{i}$. We have proved the :
\begin{lem}
\label{lem54} Assume that $\alpha_i$ and $t_i$ satisfy
(\ref{eq5.001}). Then, there exists $W_i$ a positive solution of
(\ref{eq5.3}) in $(T_{2i},T_{2i+2})$ with boundary conditions
$W_i(T_{2i+2})=v_{p,i}(T_{2i+2}+t_i)$ and $\partial_t
W_i(T_{2i+2})= \partial_t  v_{p,i}(T_{2i+2}+t_i)$. In addition, we
have the estimates
$$
\|W_i - v_{p,i}(\cdot+t_i)\|_{{\cal C}^0_\delta } \leq c \,
\lambda \, e^{(-\delta-2)T_{2i+1}}
$$
where the constant $c$ is independent of $\varepsilon$ and of the
parameters $\alpha_i$, $t_i$ and $\xi$.
\end{lem}
Observe that the solution we have obtained is unique and depends
continuously on the parameters $\alpha_i$, $t_i$ and $\xi$ since
it is the unique solution of an ordinary differential equation.
This fact is even true when $i=\ell-1$ even though the solution is
defined on a half line.

\medskip

\noindent {\bf Step 2.} We now look for a solution of
(\ref{eq5.3}) which is defined on $(T_0,T_2)$.  We decompose
$$
W_0 (t)=v_{p,0}(t+t_0)+ \overline{w}(t)+\underline{w}(t),
$$
where $\overline w$ is the solution of
\[
\label{eq5.16} L_{p,0} \, \overline{w}=-\lambda
e^{-2t}v_{p,0}(t+t_0)
\]
in $(T_0,T_2)$ with boundary data $ \overline{w}(T_2)=\partial_t
{\overline{w}}(T_2)=0$. The operator $L_{p,0}$ is the one which
has been defined in (\ref{eq5.4}). With this in mind, it remains
to find a  $\underline w$ solution of
\[
\label{eq5.17} L_{p,0}\underline{w} = - \lambda
e^{-2t}(\overline{w}+\underline{w}) - Q_0(\underline{w})
\]
in $(T_0,T_2)$ with  boundary data $\underline{w}(T_2)=\partial_t
{\underline{w}}(T_2)=0$, where
$$
Q_0(\underline w) : = |v_{p,0}(\cdot +t_0) + \overline{w} +
\underline w |^{p-1}(v_{p,0}(\cdot +t_0) + \overline{w} +
\underline w ) - v_{p,0}^p ( \cdot+t_0) - p \, v_{p,0}^{p-1}(
\cdot + t_0)(\underline w + \overline{w} ).
$$
It will be convenient to define
\[
q : = \min \left\{{2N+\delta-6\over N-2},{(N+\delta-2)p\over N-2}
\right\}.
\]
Observe that we have $q > {1\over 2}$ since we have assumed that
$N\geq 5$ and $\delta \in (-{N-1\over 2},-{N-2\over 2})$. This
time we consider the following set of functions
$$
{\cal E}_{\kappa ,0}= \left\{ w\in {\cal C}^0_\delta ((T_0, T_2))
\quad : \quad  \|w\|_{{\cal C}^0_\delta} \leq  \kappa \,
\varepsilon^q \right\} ,
$$
where the constant $\kappa >0$ will be fixed later on. It is clear
that
\[
\label{eq5.18} \|\lambda e^{-2t}v_{p,0}\|_{{\cal C}^0_\delta} \leq
c \, \lambda \, e^{(-\delta-2) \, T_1} \leq c \, \lambda \,
\varepsilon^{2 + \delta \over N-2}
\]
Using the result of Proposition \ref{lem5.1}, we get
\[
\label{eq5.19} \|\overline{w}\|_{{\cal C}^0_\delta} \leq c \,
\|\lambda e^{- 2 t} \, v_{p,0} \|_{{\cal C}^0_\delta} \leq c \,
\mu \, \varepsilon^{N-2+\delta\over N-2}\leq c \,
\varepsilon^{N-2+\delta\over N-2}
\]
so that \beq \label{eq5.20} \|\lambda e^{-2t}\overline{w}\|_{{\cal
C}^0_\delta}\leq c \, \varepsilon^{2N-6+\delta\over N-2}. \eeq As
in Step 1, we have \beq \label{eq5.21} \|\lambda
e^{-2t}{w}\|_{{\cal C}^0_\delta}\leq \lambda \|{w}\|_{{\cal
C}^0_\delta}\leq c \, \kappa \, \varepsilon^{{N-4\over N-2} + q}.
\eeq for all ${w}\in {\cal E}_{\kappa ,0}$. For $\varepsilon$
small enough, we have $1 < p < 2 $. Thus, for all $s_2\in\R$ and
all $s_1>0$, we can write
$$
\left| |s_1+s_2|^{p-1}(s_1+s_2) -s_1^p - p \, s_1^{p-1}s_2 \right| \leq c \,
|s_2|^p.
$$
for some constant $c >0$. Consequently, we can estimate for all $w
\in {\cal E}_{\kappa,0}$ \beq \label{eq5.22}
\begin{array}{rlllll}
\|Q_0(w)\|_{C^0_\delta} & \leq & c \, \sup_{t\in (T_0,T_2)} \,
e^{-\delta t} \, | w(t)+ \overline{w}(t)|^p \\[3mm]
& \leq &  c \, ( \kappa^p \, \varepsilon^{p \, q} +
\varepsilon^{p(N-2+\delta)\over N-2} ).
\end{array}
\eeq
Using the result of Proposition~\ref{lem5.1}, we get a solution
$\tilde w$ of
\[
L_{p,0} \tilde w  = - \lambda \, e^{-2t} \, (\overline{w} + w) -
Q_0(w)
\]
with $ \tilde w (T_2) = \partial_t \tilde w (T_2)=0$. Collecting
(\ref{eq5.20}), (\ref{eq5.21}) and (\ref{eq5.22}) we get the
estimate \beq \label{eq5.23} \| \tilde w  \|_{C^0_\delta}\leq
\varepsilon^{q} \, \tilde c \, ( 1 + \kappa^p \,
\varepsilon^{(p-1) \, q} + \kappa \, \varepsilon^{N-4\over N-2}).
\eeq We choose the constant $\kappa$ so that
$$
\tilde c \, (1+ \kappa^p \, \varepsilon^{(p-1) \, q}+\kappa \,
\varepsilon^{N-4\over N-2})\leq \kappa ,
$$
for all $\varepsilon $ is close to $0$, say $\varepsilon \in (0,
\varepsilon_0)$.

\medskip

As Step 1, we can define the mapping
$$
{\cal T}_{0}: {\cal E}_{\kappa,0} \rightarrow {\cal E}_{\kappa,0}
$$
by ${\cal T}_{0}(w):= \tilde w$. Clearly, ${\cal T}_{0}$ is well
defined and is  continuous and compact, so that one can again
refer to Schauder's fixed point Theorem to obtain the fixed point
of ${\cal T}_0$. We have proved the~:
\begin{lem}
\label{lem55} Given $\alpha_0$ and $t_0$ satisfying
(\ref{eq5.001}), there exists a solution $W_0$ of (\ref{eq5.3}) in
$(T_{0},T_{2})$ with  boundary conditions $W_0(T_{2}) = v_{p,0}
(T_{2} + t_0)$ and $\partial_t W_0(T_{2}) = \partial_t v_{p,0}
(T_{2} + t_0)$. In addition, we have the estimates
$$
\|W_0  - v_{p,0}(\cdot+t_0)-\bar w  \|_{{\cal C}^0_\delta} \leq  c
\, \varepsilon^q
$$
for some $q > {1\over 2}$ and
$$
\|\overline{w}\|_{{\cal C}^0_\delta} \leq c \,
\varepsilon^{N-2+\delta\over N-2}
$$
\end{lem}
Again this solution is unique and depends continuously on the
parameters $\alpha_0$, $t_0$ and $\xi$.

\medskip

\noindent {\bf Step 3.} We now explain how to choose the
parameters $(\alpha_0, \ldots, \alpha_{\ell-2})$ and $(t_0, \ldots
,t_{\ell-2})$ so that the Cauchy data of $W_i$ and $W_{i-1}$
coincide at $T_{2i}$, for $1\leq i\leq \ell-1$. To this aim, we
argue inductively, starting by matching the Cauchy data of
$W_{\ell-2}$ and $W_{\ell-1}$.

\medskip

This amounts to find $\alpha_{\ell-2}$ and $t_{\ell-2}$ so that
\[
\label{eq5.32} W_{\ell-2}(T_{2\ell-2})=
W_{\ell-1}(T_{2\ell-2}), \qquad \mbox{and} \qquad
\partial_t W_{\ell-2}(T_{2\ell-2})=\partial_t  W_{\ell-1}(T_{2\ell-2}),
\]
In other words, we need to find $\alpha_{\ell-2}$ and $t_{\ell-2}$
so that \beq \label{eq5.33}
v_{p,\ell-2}(T_{2\ell-2}+t_{\ell-2})=W_{\ell-1}(T_{2\ell-2})\qquad
\qquad \partial_t v_{p,\ell-2}(T_{2\ell-2}+t_{\ell-2}) =
\partial_t W_{\ell-1}(T_{2\ell-2}),
\eeq
It follows from Lemma \ref{lem2.5} and Lemma \ref{lem54} that
\[
\begin{array}{rlllll}
W_{\ell-1}(T_{2\ell-2}) & = & \varepsilon_{p,1} +
w_{\ell-1}(T_{2\ell-2}) & = & \varepsilon_{p,1} + {\cal O}
(\varepsilon^{\delta+N+4\ell-6\over
N-2})\\[3mm]
\partial_t  W_{\ell-1}(T_{2\ell-2}) & = & \partial_t
w_{\ell-1}(T_{2\ell-2}) & =  & {\cal O}
 (\varepsilon^{\delta+N+4\ell- 6 \over N-2}).
\end{array}
\]
and we also have
\beq \label{eq5.35}
\begin{array}{rllllll}
W_{\ell-2}(T_{2\ell-2}) & = & (\varepsilon_{p,1}
+\alpha_{\ell-2})\cosh
({N-2\over 2}t_{\ell-2}) + F_{\ell-2}\\[3mm]
\partial_t W_{\ell-2} (T_{2\ell-2}) & = &
(\varepsilon_{p,1}+\alpha_{\ell-2}){N-2\over 2} \sinh ({N-2\over
2}t_{\ell-2}) + G_{\ell-2} ,
\end{array}
\eeq where the continuous functions $(F_{\ell-2},G_{\ell-2})$
depend on $\alpha_{\ell-2}$ and $t_{\ell-2}$ and satisfy
$|F_{\ell-2}|+|G_{\ell-2}|={\cal O}( \varepsilon^{p\over 2})$. The
system (\ref{eq5.33}) is therefore equivalent to \beq
\label{eq5.36} t_{l-2} = \varepsilon^{-{1\over 2}} \, \tilde
F_{l-2}\qquad \mbox{and} \qquad \alpha_{l-2} = \tilde G_{l-2} \eeq
where the continuous functions $(\tilde F_{l-2}, \tilde G_{l-2})$
depend on $\alpha_{\ell-2}$ and $t_{\ell-2}$ and satisfy $|\tilde
F_{\ell-2}|+| \tilde G_{\ell-2}|={\cal O}( \varepsilon^{p\over
2})$ (here we have used the fact that ${\delta+N+4l-6\over
N-2}>{p\over 2}$ provided $ \varepsilon $ is close to $0$).

\medskip

Recall that $p=p_N+\varepsilon$. Given $\gamma \in(1, {p\over 2})$ we
define
$$
B := \left\{ (\alpha_{\ell-2}, t_{l-2}) \in \R^2 \quad  :  \quad
\alpha_{l-2}^2 + \varepsilon \, t_{l-2}^2 \leq \varepsilon^\gamma
\right\}.
$$
In view of (\ref{eq5.36}), the mapping
\[
\tilde H : B \longrightarrow B
\]
defined by $ \tilde H (\alpha_{\ell-2}, t_{\ell-2})  : =
(\varepsilon^{-{1\over 2}} \, \tilde F_{\ell-2}, \tilde G_{\ell-2}
)$ is well defined and it follows from Browder's fixed point
theorem that  (\ref{eq5.36}) admits a solution. In addition,
applying Lemma~\ref{lem2.5} and Lemma~\ref{lem2.6}, we get \beq
\label{eq5.37}
v_{p,\ell-2}(T_{2\ell-4}+t_{\ell-2})=\varepsilon_{p,2}+ {\cal O} (
\varepsilon^{p\over 2}), \qquad \qquad \partial_t
v_{p,\ell-2}(T_{2\ell-4}+t_{\ell-2})= {\cal O} (
\varepsilon^{p\over 2}). \eeq

\medskip

Arguing inductively, we construct a function $v_{p,\lambda,\xi}$,
solution of (\ref{eq5.3}), which depends on $\lambda$ and on
$\xi$, and which satisfies \beq \label{eq5.38} v_{p,\lambda,\xi}
(T_2)=\varepsilon_{p,\ell-1}+ {\cal O} ( \varepsilon^{p\over 2}),
\qquad \qquad
\partial_t v_{p,\lambda,\xi}(T_2)= {\cal O} ( \varepsilon^{p\over 2}).
\eeq
In view of (\ref{eq2.96}) and thanks to Lemma 1, Lemma 7, Corollary 3 and Corollary 4,
the following expansion holds
\beq
\begin{array}{ll}
&v_{p, \lambda, \xi} (t)\\
=&\ds{(\ell\varepsilon)^{1\over 2}\left[{\sqrt{ C_N^{(4)}}\over 2}e^{(N-2)(\xi-t)\over 2}
+{\sqrt{ C_N^{(4)}}\over 2}e^{(2-N)(\xi-t)\over 2}-{4\mu  C_N^{(8)}\ell^{4-N\over N-2}\over
(N-2)^2\sqrt{ C_N^{(4)}}}e^{(N-6)\xi\over 2}e^{(2-N)t\over 2}
\right]}\\
&+{\cal O}\left(\varepsilon^q e^{\delta t}+\varepsilon^{{1\over 2}+{3\over N+2}}+
(\varepsilon^{1\over 2}e^{(N-2)t\over 2})^p+ \varepsilon^{1\over 2}(\varepsilon^{2\over N}\log{1\over\varepsilon }
+\varepsilon^{N-17/4\over N-2}) e^{(2-N)t\over 2}\right)
\end{array}
\eeq
for all $t\in (\log {1\over 2r_\varepsilon},\log {2\over r_\varepsilon})$.
This completes the proof of the result. \cqfd

\medskip

\section{The linear analysis.}

Assume that $\O$ is a regular bounded open subset of $\R^N$ and
$\Sigma:=\{a_1,\ldots , a_m\}$ is a finite set of points of $\O$.
We choose $r_0>0$ in such a way that the closed balls $B(a_i , 2
r_0)$, for $i=1, \ldots, m$ are disjoint and included in $\O$. For
all $r \in (0, r_0)$, we define
\[
\O_{int,r} : = \bigcup_{i=1}^m B(a_i,r) \qquad \mbox{and} \qquad
\O_{ext,r} := \O -  \overline \O_{int,r} .
\]
We define the weighted spaces~:
\begin{defi}
Given $\nu \in \R$, the space ${\cal C}^0_\nu (\overline{\O} -
\Sigma)$ is defined to be the set of continuous functions $w \in
C^{0}_{loc} (\overline \O - \Sigma )$ for which the following norm
is finite~:
\begin{equation}
\label{eq3.0} \| w \|_{C^0_\nu (\overline \O - \Sigma )} : =  \|
w\|_{{L^\infty}(\overline{\O_{ext,r_0}})} + \sum_{j=1}^m \|
r^{-\nu} \, w ( a_j + \cdot) \|_{L^\infty (B(0,2r_0) -\{0\})}.
\end{equation}
\end{defi}
Given $r \in (0, r_0)$ we define the space $C^0_\nu (\overline
\Omega_{ext,r})$ to be the space of restrictions of functions of
$C^0_\nu (\overline{\O} - \Sigma)$ to $\overline \Omega_{ext,r}$.
This space is endowed with the induced norm.

\medskip

In this section, we study the linearization of the nonlinear
operator (\ref{eq5.2}) about the radial function
\[
u_{\e} (x) := |x|^{-{2\over p-1 }}   \, v_{\e} (-\log |x|)
\]
where $v_\e : = v_{p,\lambda,\xi}$ and $v_{p,\lambda,\xi}$
is the solution of (\ref{eq5.3}) defined in Step 3 of the proof of
Proposition 7. This operator is defined by
\[
L_\e := \Delta  + \lambda + p \, u_{\e}^{p-1}.
\]
Recall $r_\varepsilon= \varepsilon^{2\over N^2-4}$. We can write any function $v$ defined in the punctured ball
$B(0,r_\varepsilon)-\{0\}$ as
\[
v(x) = |x|^{-\frac{2}{p-1}} \, w(-\log |x|,\theta),
\]
so that the study of $L_\e$ reduces to the study of the linear
operator
\begin{equation}
{\cal L}_\e := \partial^2_t - a_p \,  \partial_t  - b_p \, +
\Delta_{S^{N-1}}  + p \,v_{\e}^{p-1}  + \lambda \, e^{-2t}
\end{equation}
on the half cylinder $[B_\varepsilon,+\infty)\times S^{N-1}$, where
$\Delta_{S^{N-1}}$ denotes the Laplace-Beltrami operator on the
sphere $S^{N-1}$ and $B_\varepsilon=-\log r_\varepsilon$.

\medskip

We denote by $ (e_j, \lambda_j)$ the set of eigendata of
$\Delta_{S^{N-1}}$
\[
\Delta_{S^{N-1}} e_j =  - \lambda_j \, e_j.
\]
We also assume that the eigenvalues are counted with multiplicity,
that $\lambda_j \leq \lambda_{j+1}$ and that the  $e_j$ are
normalized by
\[
\int_{S^{N-1}} e_j^2 \, d\omega = 1.
\]

We now prove some uniform estimates for a right inverse for the
operator ${\cal L}_\e $.
\begin{pro}
\label{pro3.1} Assume that $ \delta \in (-{N+2\over 2},- {N\over
2})$ is fixed. Then, there exists $p_0 \in (p_N,+\infty)$ such
that, if $p \in ( p_N,p_0)$, then, for all $f\in {\cal C}^0_\delta
([B_\varepsilon,+\infty)\times S^{N-1})$, there exists a unique solution $w
\in {\cal C}^0_\delta ([B_\varepsilon,+\infty)\times S^{N-1})$ of \beq
\label{eq4.0} {\cal L}_{\e} \, w = f \eeq in $[B_\varepsilon,+\infty)\times
S^{N-1}$ which satisfies
\[
w(B_\varepsilon,\theta) \in \mbox{Span} \{ e_j \, : \, j=0, \ldots, N\}.
\]
Furthermore, \beq \label{eq4.01} \|w\|_{{\cal C}^0_\delta} \leq c
\, \|f\|_{C^0_\delta} \eeq for some constant which does not depend
on $\e$.
\end{pro}
{\it Proof.} The proof is divided in three parts. In the first
part we explain how to solve the equation (\ref{eq4.0}) when the
function $f$ does not have any component on $e_j$ for $j=0,
\ldots, N$ in its eigenfunction decomposition. Next, in the second
part, we obtain a uniform estimate for the solution already
obtained. Finally, in the last part, we explain how to solve
(\ref{eq4.0}) when the eigenfunction decomposition of $f$ has
components on $e_0, \ldots, e_N$.

\medskip

\noindent {\bf Step 1} For the time being, we assume that the
eigenfunction decomposition of the function $f$ is given by \beq
\label{eq4.1} f(t,\theta)=\sum_{j\geq N+1} f^j(t) \, e_j(\theta).
\eeq

Observe that, as $p$ tends to $p_N$ we have
\[
\lim_{p\rightarrow p_N} a_p= 0 ,\qquad \qquad \lim_{p\rightarrow
p_N} b_p  = \left( \frac{N-2}{2} \right)^2
\]
and
\[
\lim_{p\rightarrow p_N} \sup \, p\, v_{p, \lambda,\xi}^{p-1} =
{N(N+2)\over 4}
\]
these limits being independent of the solution $v_{p, \lambda,\xi}$.

\medskip

Now the eigenfunction decomposition of the Laplace-Betrami
operator on $S^{N-1}$ induces a decomposition of the operator
${\cal L}_{\e}$  into the sequence of operators
\[
L_{j} := \partial^2_t - a_p \,  \partial_t  - b_p \,
 - \lambda_j + p \,v_\e^{p-1}  + \lambda \, e^{-2t} .
\]

Using these above limits together with the fact that $\lambda_j
\geq 2\, N$ for $j\geq N+1$, we conclude that, for $j \geq N+1$
the potential is negative provided $p$ is close enough to $p_N$.
In particular, this implies that it is possible to solve
\[
{\cal L}_{\e} \, w =f
\]
on any $(B_\varepsilon, S)\times S^{N-1}$, with $w=0$ as boundary data
(observe that the operator ${\cal L}_{\e}$ is not self adjoint but
is conjugate to a self adjoint operator and we have just seen that
this former is injective, when restricted to the set of functions
spanned by $e_j$, for $j \geq N+1$).

\medskip

It remains to prove that there exists a constant $c >0$ which does
not depend on $S$, nor on $p$ such that
\begin{equation}
\label{eq:qsqs}
\sup |e^{-\delta t} \, w| \leq c \, \sup |
e^{-\delta t} \, f|.
\end{equation}
Then, the existence of the solution on all $(B_\varepsilon, +\infty) \times
S^{N-1}$ as well as the relevant estimate will follow by passing
to the limit $S\rightarrow +\infty$. To keep the proof short and
since anyway our aim is to pass to the limit as $S$ tends to
$\infty$, it is enough to prove that (\ref{eq4.0}) holds for all
$S$ chosen large enough so that $\sup_{(S, +\infty)} v_p \leq \e$.

\medskip

\noindent {\bf Step 2} The proof of (\ref{eq:qsqs}) is by
contradiction. If it were false for all choice of $p_0$ and $S$,
there would exist a sequence $(p_n)_n$ tending to $p_N$, a
sequence of functions  $(f_n)$ and a sequence of reals $(S_n)_n$
and a sequence $(w_n)_n$ of solutions of (\ref{eq4.0}) such that
\beqr \|f_n\|_{{\cal C}^0_\delta } \equiv 1
 \qquad \mbox{and}\qquad
\lim_{n\rightarrow +\infty} A_n: =  \sup e^{-\delta t} |w_n| =
+\infty. \eeqr
We denote $B_n=B_{\e_n}$ where $\e_n:=p_n-p_N$. Obviously, there exists a point $(t_n,\theta_n)\in
(B_n,S_n)\times S^{N-1}$ where the above supremum is achieved,
namely $A_n = e^{- \delta t_n} \, |w_n(t_n,\theta_n)|$. Observe
that elliptic estimates imply that \begin{equation} \sup
e^{-\delta t} |\nabla w_n| \leq c \, (1+ A_n) \end{equation} and
this in turn implies that the sequences $(t_n-B_n)_n$  and
$(S_n-t_n)_n$ remain bounded away from $0$.

\medskip

We define $\tilde t_n >B_n$ to be the nearest local maximal point of
the function $v_{p_n}(t)$ to the point $t_n$. We distinguish
several cases according to the behavior of the sequence
$(t_n)_n$.

\medskip

\noindent{\bf Case 1.} Assume that the sequence $(t_n-\tilde
t_n)_n$ is bounded. In this case, we define the function $\tilde
w_n$ by
\[
\tilde w_n (t,\theta)= \frac{1}{A_n} \, e^{- \delta \tilde t_n}\,
w_n (t + \tilde t_n, \theta).
\]
Observe that the sequence of functions $(v_{p_n} (\cdot + \tilde
t_n))_n$ converges on compact to $t\rightarrow
(N(N-2))^{\frac{N-2}{4}}(\cosh t)^{\frac{2-N}{2}}$ (see
Proposition 5). Up to a subsequence, we can assume that the
sequence $(t_n -\tilde t_n)_n$ converges to $t_\infty$. Moreover,
we can assume that the sequence $(\tilde w_n)_n$ converges on
compacts to $\tilde w_\infty$ a nontrivial solution of \beq
\label{eq4.002}
\partial^2_t \tilde w_\infty + \Delta_{S^{N-1}} \, \tilde
w_\infty - {(N-2)^2\over 4} \, \tilde w_\infty + \frac{N(N+2)}{4}
\, (\cosh t)^{-2} \, \tilde w_\infty = 0.
\eeq
Moreover, $\tilde w_\infty$ is bounded by a constant times
$e^{\delta t}$. The fact that $\tilde w_\infty$ is not identically
equal to $0$ follows from the fact that $|\tilde w_n (t_n -\tilde
t_n, \theta_n)| = e^{\delta (t_n-\tilde t_n)}$ and hence remains
 bounded away from $0$.

\medskip

We consider the eigenfunction decomposition of $\tilde w_\infty$
\[
\tilde w_\infty = \sum_{j=N+1}^\infty a_j \, e_j .
\]
At $-\infty$ the function $a_j$ is either blowing up like $t
\longrightarrow e^{- \gamma_j t}$ or decaying like $t
\longrightarrow e^{\gamma_j t}$, where
\[
\gamma_j : = \sqrt{\lambda_j+{(N-2)^2\over 4}} .
\]
The choice of $\delta \in (-\frac{N+2}{2}, - \frac{N}{2})$ implies
that $- \delta < \gamma_j$ for all $j \geq N+1$. Hence $a_j$
decays exponentially at $-\infty$. Multiplying the equation (\ref{eq4.002})
by $a^j \, e_j$ and integrating by parts over $\R$ (all
integrations are justified because $a_j$ decays exponentially at
both $\pm \infty$), we get
$$
\begin{array}{rllll}
\ds{\int_{-\infty}^{+\infty}|\partial_t a^j|^2 + (\lambda_j +
{(N-2)^2\over 4}) (a^j)^2} & = & \ds \frac{N(N+2)}{4} \,
\int_{-\infty}^{+\infty} (\cosh s)^{-2} \, (a^j)^2 \\[3mm]
&  \leq & \ds \frac{N(N+2)}{4} \int_{-\infty}^{+\infty}(a^j)^2.
\end{array}
$$
Since $j\geq N+1$, we have  $\lambda_j\geq 2N$, and hence we
conclude that $a_j\equiv 0$. Hence, $\tilde w_\infty\equiv 0$, a
contradiction.

\medskip

\noindent{\bf Case 2.} Assume that the sequence $(t_n-\tilde
t_n)$,  the sequence $(t_n-B_n)_n$  and the sequence $(S_n-t_n)_n$ are
all unbounded. In this case, we define the function $\tilde w_n$
by
\[
\tilde w_n (t,\theta)= \frac{1}{A_n} \, e^{- \delta t_n}\, w_n (t
+ t_n, \theta).
\]
Observe that this time the sequence of functions $(v_{p_n} (\cdot
+ t_n))_n$ converge to $0$ on compacts. Up to a subsequence, we
can assume that the sequence $(\tilde w_n)_n$ converges on
compacts to $\tilde w_\infty$ a nontrivial solution of
\[
\partial^2_t \tilde w_\infty + \Delta_{S^{N-1}} \, \tilde
w_\infty - {(N-2)^2\over 4} \, \tilde w_\infty = 0.
\]
Moreover, $\tilde w_\infty$ is bounded by a constant times
$e^{\delta t}$.

\medskip

Again, we consider the eigenfunction decomposition of $\tilde
w_\infty$
\[
\tilde w_\infty = \sum_{j=N+1}^\infty a_j \, e_j
\]
and we see that $a_j$ is a linear combination of $t
\longrightarrow e^{- \gamma_j t}$ and $t \longrightarrow
e^{\gamma_j t}$. The choice of $\delta \in (- \frac{N+2}{2}, -
\frac{N}{2})$ implies that $\delta > - \gamma_j$ for all $j \geq
N+1$. Hence $a_j$ cannot be bounded by $e^{\delta t}$ unless it is
identically $0$. We conclude that $a_j\equiv 0$. Hence, $\tilde
w_\infty\equiv 0$, a contradiction.

\medskip

\noindent {\bf Case 3.} Assume that the sequence $(t_n-B_n)_n$ is
bounded (resp. that the sequence $(S_n-t_n)_n$ is bounded) and
that the sequence $(t_n-\tilde t_n)$ is unbounded. This case can
be treated as in case 2. The only difference is that this time
$\tilde w_\infty$ is defined on $[\underline t_\infty, +\infty)
\times S^{N-1}$  (resp. on $(-\infty, \bar t_\infty]\times
S^{N-1}$) and is equal to $0$ on $\{\underline t_\infty\} \times
S^{N-1}$ (resp. on
$
\{\bar t_\infty\}\times S^{N-1}$). We omit the details.

\medskip

Since we have reached a contradiction in each case, the proof of
the claim is complete. We can now pass to the limit as $S$ tends
to $+\infty$ and complete the proof of the result in the case
where the eigenfunction decomposition of $f$ does not involve any
$e_j$ for $j=0, \ldots, N$.

\medskip

\noindent {\bf  Step 3.} Now we consider the case where the
function $f$ is collinear to $e_j$, namely
\[
f(t,\theta)=f^j(t)\, e_j(\theta)
\]
for some $0\leq j\leq N$. We extend the function $f$ to
be equal to $0$ when $t\leq B_\e$ and we define the function $\tilde
v_p$ which is equal to $v_p$ for $t\geq B_\e$ and is equal to $0$
for $t <B_\e$. We consider the equation
\beq
\label{eq4.200} \partial_t^2 a^j - a_p \,
\partial_ta^j-(\lambda_j+b_p) \, a^j + p \, \tilde v_p^{p-1} \,
a^j+\lambda e^{-2t} \, \chi \, a^j = f^j \eeq in ${\R}$. Here
$\chi $ is a cutoff function identically equal to $1$ on
$(B_\e,+\infty)$ and equal to $0$ on $(-\infty, B_\e-1)$. Observe that
$$
|f^j_n(t)| \leq \|f \|_{{\cal C}^0_\delta} \, e^{\delta t}.
$$
For $p$ close enough to $p_N$, $\delta$ is not an indicial root of
the operator ${\cal L}_\e$ and it follows from Cauchy's theorem
that there exists a unique solution of (\ref{eq4.200}) which is
bounded by a constant times $e^{\delta t}$ at $+\infty$. A priori
this solution is only defined for $t$ large enough but is can be
extended to all $\R$ easily. Furthermore, it follows from the
construction of this solution that
\[
\sup_{(T, +\infty)} e^{- \delta t} \, |a^j| \leq c\, \sup_{\R}
e^{- \delta t} \, |f^j|
\]
provided $T$ is large enough. This solution satisfies
\beq
\label{eq4.201}
 \partial_t^2 a^j - a_p \,
\partial_t a^j-(\lambda_j+b_p) \, a^j =0
 \eeq
for $t <B_\e-1$ and, since $\delta \in (-\frac{N+2}{2},-\frac{N}{2})$, even if
$a^j$ blows up at $-\infty$, it blows up at a slower rate than  $t
\rightarrow e^{\delta t}$, provided $p$ is chosen close enough to
$p_N$.

\medskip

We claim that there exists a constant $c >0$ such that
\[
\sup_{\R} e^{- \delta t} \, |a^j| \leq c \, \sup_{\R} e^{- \delta
t} \, |f^j|
\]
provided $p$ is close enough to $p_N$. As before, we argue by
contradiction. Assume that the claim is not true. Then there would
exist a sequence $(p_n)_n$ tending to $p_N$, a sequence of
functions $(f^j_n)_n$ and a sequence of solutions $(a^j_n)_n$ of
(\ref{eq4.201}) such that
\[
\sup_{\R} e^{- \delta t} \, |f^j_n|  =1 \qquad \mbox{and} \qquad
A_n : = \sup_{\R} e^{- \delta t} \, |a^j_n|
\]
tends to $+\infty$. The previous remarks show that the above
supremum is always achieved in ${\R}$. So we can define $t_n$ such
that $A_n = e^{- \delta t_n} \, |a^j_n(t_n)| $.

\medskip

As in Step 2, we define $\tilde t_n >0$ to be the nearest local
maximal point of the function $v_{p_n}(t)$ to the point $t_n$. We
distinguish several cases according to the behavior of the
sequence $(t_n)_n$. We define the function $\tilde a_n^j$ by
\[
\tilde a^j_n (t)= \frac{1}{A_n} \, e^{- \delta t_n}\, a_n^j (t +
\tilde t_n) .
\]

We can assume that, up to a subsequence, the sequence
$(\tilde a^j_n)_n$ converges on compacts to $\tilde a_\infty$ a
nontrivial solution of
\[
\partial^2_t \tilde a_\infty - \lambda_j \, \tilde
a_\infty - {(N-2)^2\over 4} \, \tilde a_\infty + \frac{N(N+2)}{4}
\, (\cosh t)^{-2} \, \tilde a_\infty = 0
\]
in the case where  the sequence $(t_n-\tilde t_n)$ is bounded, or
to a nontrivial solution of
\[
\partial^2_t \tilde a_\infty -\lambda_j \, \tilde
a_\infty - {(N-2)^2\over 4} \, \tilde a_\infty = 0
\]
in the case where  the sequence $(t_n-\tilde t_n)$ is unbounded.

\medskip

Moreover, $\tilde a_\infty$ is bounded by a constant times
$e^{\delta t}$. However, the choice of $\delta \in
(-\frac{N+2}{2}, - \frac{N}{2})$ implies that $\delta < -\gamma_j$
for all $j =0, \ldots ,N$ and there are non nontrivial solutions
of the above homogeneous problems which are bounded by $e^{\delta
t}$ at $+\infty$. Hence, $\tilde a_\infty\equiv 0$, a
contradiction.  This completes the proof of the result. \cqfd

\medskip

We recall some well known result concerning harmonic extension of
functions which are defined on $S^{N-1}$.
\begin{lem}
Given $\varphi\in C^{2,\alpha}(S^{N-1})$, we define $V_\varphi$ to
be the unique harmonic extension of $\varphi$ in $B(0,1)$, namely
\begin{equation}
\left\{
\begin{array}{rllll}
\Delta V_\varphi & = & 0 \qquad & \mbox{in } \qquad
B(0,1)\\[3mm]
V_\varphi& = & \varphi \qquad & \mbox{on } \qquad \partial B(0,1)
\end{array}
\right.
\end{equation}
Assume that $\varphi$ is $L^2(S^{N-1})$ orthogonal to $e_0,
\ldots, e_N$, then
\[
\| V_\varphi \|_{{\cal C}^{0}_2 (B(0, 1)-\{0\})} \leq c \,
\|\varphi\|_{{\cal C}^{0}(S^{N-1})}
\]
for some constant $c >0$ which does not depend on $\varphi$.
\label{lem6.1}
\end{lem}

Using the fact that Kelvin's transform of an harmonic function $V$
\[
W(x) = |x|^{2-N} \, V(\frac{x}{|x|})
\]
is harmonic, the above result translates into the~:
\begin{lem}
\label{lem6.2} Given $\varphi\in C^{2,\alpha}(S^{N-1})$, we define
$W_\varphi$ to be the unique harmonic extension of $\varphi$ in
${\R}^N-B(0,1)$ which decays at $\infty$. Namely
\begin{equation}
 \left\{
\begin{array}{rllll}
\Delta W_\varphi & = & 0 \qquad &\mbox{in } \qquad \R^N-
B(0,1)\\[3mm]
W_\varphi & = & \varphi \qquad &\mbox{on } \qquad \partial B(0,1)
\end{array}
\right.
\end{equation}
and $W_\varphi$ tends to $0$ at $\infty$. Assume that $\varphi$ is
$L^{2}(S^{N-1})$ orthogonal to $e_0,\cdots,e_N$ then
\[
\| W_\varphi \|_{{\cal C}^{0}_{-N} ({\R}^N-B(0,1))} \leq c
\, \|\varphi\|_{{\cal C}^{0} (S^{N-1})}
\]
for some constant $c>0$ which does not depend on $\varphi$.
\end{lem}
From now on we assume that $\O$ is a bounded regular domain in
$\R^N$.

\section{Bubble tree solutions in general domains}

As before, we only prove the case when $N\geq 6$ since the proof
of the result when $N=5$ follows the same lines with minor
modifications. We recall
\[
r_\varepsilon = \varepsilon^{2\over N^2-4}.
\]

We define the space
\[
{\cal E} : = \left\{ \varphi \in C^{2,\alpha}(S^{N-1}) \quad :
\quad \int_{S^{N-1}} \varphi  \, e_j \, d\omega= 0, \quad  j= 0 ,
\ldots , N \quad \mbox{and} \quad  \|\varphi\|_{{\cal C}^{2,
\alpha}} \leq r_\varepsilon \, \varepsilon^{1 \over 2} \right\}.
\]

\subsection{Solution of the nonlinear problem in
$\Omega_{int,\varepsilon}$.}

Given a $m$ functions ${\bf\varphi} := (\varphi_1,
\ldots,\varphi_m) \in {\cal E}^m$ and $m$ points ${\bf x} : =
(x_1, \ldots, x_m) \in \O^m$, we construct a positive solution of problem
(\ref{eq5.2}) in $\Omega_{int,\varepsilon}$ whose boundary is, in
some sense, parameterized by ${\bf \varphi}$. Namely we would like
to solve
\begin{equation}
\left\{
\begin{array}{rlllll}
\Delta  u_{int,i} + \lambda \, u_{int,i} + u_{int,i}^p & = & 0
\qquad &\mbox{in } \qquad
B(x_i,r_\varepsilon)\\[3mm]
u_{int,i} & \in  & \mbox{Span} \, \{ e_0, \ldots, e_N\} \qquad
&\mbox{on } \qquad \partial B(x_i,r_\varepsilon)
\end{array}
\right. \label{eq6.3.00}
\end{equation}

\medskip

For each $i=1, \ldots, m$, we denote by $V_{\varphi_i}$ the unique
harmonic extension of $\varphi_i$ in $B(x_i,R_\varepsilon)$,
namely
\begin{equation}
\left\{
\begin{array}{rlllll}
\Delta  V_{\varphi_i} & = & 0 \qquad &\mbox{in } \qquad
B(x_i,r_\varepsilon)\\[3mm]
V_{\varphi_i} & = & \varphi_i \qquad &\mbox{on } \qquad \partial
B(x_i,r_\varepsilon)
\end{array}
\right.
\label{eq6.3}
\end{equation}
It follows from Lemma~\ref{lem6.1}, together with a scaling
argument, that
\begin{equation}
\| V_{\varphi_i} \|_{{\cal C}^{0}_2(B(x_i,
r_{\varepsilon}) -\{ x_i \}) } \leq  c \, r_\varepsilon^{-2} \, \|
\varphi_i \|_{{\cal C}^{0} (S^{N-1})}. \label{eq6.3.1}
\end{equation}

We keep the notations of the previous sections and, we look for a
positive solution of problem (\ref{eq5.2}) in $B(x_i, r_\varepsilon )$ of
the form
\begin{equation}
\label{eq6.4} u_{int,i} = u_{p,\lambda,\xi_i}(\cdot - x_i)+
V_{\varphi_i} + w_i
\end{equation}
where the function $u_{p,\lambda,\xi _i}$ is the radial solution
of problem (\ref{eq5.2}) which has been obtained in Proposition 7
and where the functions $w_i$ is small.

\medskip

As usual, we introduce the polar coordinates $(t, \theta) \in
(-\log r_\varepsilon, +\infty ) \times S^{N-1}$ in each $B(x_i,
r_\e)$. Given a function $v$, defined on $B(x_i, r_\e)$, we agree
that the function $\tilde v$ is the function defined on $(-\log
r_\varepsilon, +\infty ) \times S^{N-1}$ which is determined by
the relation
\begin{equation}
\label{eq6.5} v(x) = |x|^{-\frac{2}{p-1}} \, \tilde v( -\log |x|,
\theta).
\end{equation}
With these notations, we need to find a function $\tilde
u_{int,i}$ and $b_0, \ldots, b_N \in \R$ such that
\begin{equation}
\partial^2_t \tilde u_{int,i} - a_p \, \partial_t \tilde u_{int,i} - b_p \, \tilde u_{int,i}
+ \Delta_{S^{N-1}} \tilde u_{int,i} = - \, \lambda  \, e^{-2t} \,
\tilde u_{int,i} - \tilde u_{int,i}^p \label{eq6.8}
\end{equation}
in $[-\log r_\varepsilon,+\infty)\times S^{N-1}$ and
\[
\tilde u_{int,i}(-\log r_\varepsilon,\theta) =
r_\varepsilon^{2\over p-1}\varphi_i(\theta) + \sum_{j=0}^N b_j \,
e_j
\]
on $S^{N-1}$.

\medskip

We will obtain a solution of this equation as a fixed point for
some contraction mapping. We fix $\delta\in (-({N\over
2}+{N\over N+2}),-{N\over 2})$ such that $({2p\over p-1}+\delta-{2\over N+2}){N+2\over 2}>2$ and we define
\begin{equation}
{\cal E}_{int, \varepsilon}: = \left\{ \tilde w \in {\cal
C}^0_\delta ([-\log r_\varepsilon,+\infty)\times S^{N-1}): \quad
\|\tilde w\|_{{\cal C}^0_\delta} \leq \kappa \,
\varepsilon^{{1\over 2}+({2 p\over p-1}+\delta-{2\over
N+2}){1\over N-2}} \right\} \label{eq6.9}
\end{equation}
where the parameter $\kappa >0$ will be fixed later on.

\medskip

We write (\ref{eq6.8}) as
\begin{equation}
\label{eq6.10.1} {\cal L} \, \tilde w_i = -\lambda e^{-2t}(\tilde
w_{i}+\tilde V_{\varphi_i})-Q_{\varphi_i}(\tilde w_i)
\end{equation}
where the linear operator ${\cal L}$ is given by
\[
{\cal L} := \partial_t^2 + \Delta_{S^{N-1}} - a_p \,  \partial_t -
b_p + p \, \tilde u_{p, \lambda, \xi_i}^{p-1}
\]
and where $Q_{\varphi_i}$ collects the nonlinear terms
\[
Q_{\varphi_i}(\tilde w_i) : = (\tilde u_{p, \lambda,\xi_i} +
\tilde V_{\varphi_i} + \tilde w_i )^p- \tilde u_{p,
\lambda,\xi_i}^p - p \, \tilde u_{p, \lambda,\xi_i}^{p-1} \,
\tilde w_i .
\]
We estimate
\begin{equation}
\label{eq6.11} \| \lambda e^{-2t} \tilde V_{\varphi_i}\|_{{\cal
C}^{0}_\delta} \leq  c \, \lambda \, \varepsilon^{1\over 2} \,
r_\varepsilon^{{2\over p-1}+3+\delta}
\end{equation}
and
\begin{equation}
\|\lambda e^{-2t}\tilde w\|_{{\cal C}^0_\delta} \leq c \, \lambda
\, \kappa \, \varepsilon^{{1\over 2}+({2 p \over
p-1}+\delta-{2\over N+2}){1\over N-2}} \, r_\varepsilon^2.
\end{equation}
In view of the asymptotic expansion of $\tilde u_{p,\lambda,
\xi_i}$  we have obtained in Proposition 7, it is easy to check
that, for all $\tilde w \in {\cal E}_{int, \varepsilon}$
\[
|\tilde w | \ll \tilde u_{p, \lambda,\xi_i}
\]
in $(-\log r_\varepsilon,+\infty)\times S^{N-1}$. Moreover, it
follows from (\ref{eq6.3.1}) that
\begin{equation}
\label{eq6.6} |\tilde V_{\varphi_i} |\leq c \, r_\varepsilon^{-2}
\, \| \varphi_i \|_{L^\infty} \, e^{- {2 p\over p-1} t} \leq  c \,
 r_\varepsilon^{-1} \, \varepsilon^{1\over 2} \, e^{- {2p\over p-1} t},
\end{equation}
in  $(-\log r_\varepsilon, +\infty) \times S^{N-1}$. Hence, we
conclude that \beq \label{eq6.7} |\tilde V_{\varphi_i} |\ll \tilde
u_{p, \lambda, \xi_i} . \eeq Taylor's expansion yields
$$
(1+t)^p - 1 - p \, t \leq c \, t^2
$$
near $t=0$. This, together with the fact that $\delta < -
\frac{2}{p-1}$, implies that
\begin{equation}
\label{eq6.14} \|Q_{\varphi_i} (\tilde w) \|_{{\cal C}^0_\delta}
\leq c \, \varepsilon^{{1\over 2}+({2p \over p-1} + \delta-
{2\over N+2}) {1\over N-2}}(1+ c_{\kappa}\varepsilon^{({p-1\over
2}-{1\over N-2})}),
\end{equation} for some constant $c_{\kappa} >0$ depending on
$\kappa$. We have used the fact that ${2 p \over p-1}+\delta-
{2\over N+2} < 1$.

\medskip

Gathering the previous estimates, we conclude that \beq
\label{eq6.15} \| -\lambda e^{-2t}(\tilde w + \tilde
V_{\varphi_i})-Q_{\varphi_i}(\tilde w)\|_{_{{\cal
C}^0_\delta}}\leq c \,  (1+c_{\kappa} \, \varepsilon^{\gamma})
\varepsilon^{{1\over 2}+({2p\over p-1} + \delta-{2\over N+2})
{1\over N-2}} \eeq where $c_{\kappa}>0$ depends on $\kappa$ and
the positive number $\gamma$ is independent of $p$.

\medskip

Given $\tilde w \in {\cal E}_{int, \varepsilon}$ we use the result
of Proposition~\ref{pro3.1} to solve
\[
{\cal L} \, \tilde v = -\lambda \, e^{-2t} \, (\tilde w_{i}+\tilde
V_{\varphi_i})- Q_{\varphi_i}(\tilde w)
\]
It follows from Proprosition~\ref{pro3.1} and the above estimate
that, given $\kappa$, there exists $\e_{0} > 0 $ (depending on
$\kappa$) such that the mapping
\[
T_{i}: {\cal E}_{int, \varepsilon} \longrightarrow {\cal E}_{int,
\varepsilon}
\]
defined by $T_i (\tilde w ) = \tilde v$ is well defined, provided
$\varepsilon \in (0, \e_0)$.

\medskip

Moreover, for all $\tilde w_1, \tilde w_2\in {\cal E}_{int,
\varepsilon}$, one can check that
\beq \label{eq6.16}
\begin{array}{rlllll}
\| T_{i}(\tilde w_1) - T_{i}(\tilde w_2)\|_{{\cal C}^{0}_\delta}
&\leq & c \, \lambda \, \|\tilde w_1-\tilde w_2\|_{{{\cal
C}^{0}_\delta}} + c \, \| Q_{\varphi_i}(\tilde w_1) -
Q_{\varphi_i}(\tilde w_2)\|_{{\cal C}^{0}_\delta}\\[3mm]
& \leq & c \, (\lambda +\varepsilon^{p-1\over 2}) \, \| \tilde w_1
- \tilde w_2\|_{{{\cal C}^{0}_\delta}}.
\end{array}
\eeq Consequently, for $p$  sufficiently close to $p_N$, the
mapping $T_{i}$ is a contraction from ${\cal E}_{int,
\varepsilon}$ into itself and hence admits a unique fixed point in
this set. This yields a solution $u_{int, i}$ of (\ref{eq6.3.00}).

\medskip

If we define the function $u_{int}$ to be equal to $u_{int, i}$ on
$B(x_i, r_\e)$, we have proven the~:
\begin{pro}
\label{pro7.7} Given ${\bf x}\in\O^m$ and ${\bf\varphi} \in {\cal
E}^m$, there exists a positive solution $u_{int}$ of (\ref{eq5.2}) in
$\Omega_{int,\varepsilon}$ satisfying boundary conditions
$$
u_{int}|_{\partial B(x_i, r_\varepsilon)} -\varphi_i\in
\mbox{Span} \{ e_j \quad : \quad j=0, \ldots, N\}
$$
for all $1\leq i\leq m$. Moreover, the sequence of solutions
$u_{int}$ blows up at each $x_i$ as $p$ tends to $p_N$ in such a
way that
\[
|\nabla u_{int}|^2 \, dx  \rightharpoonup C_N^{(3)}  \sum_{i=1}^m
\ell_i \, \delta_{x_i}
\]
in the sense of measures. Here $C_N^{(3)}$ is the constant defined
in Theorem \ref{th6.1}. Finally, this solution can be expanded as
\[
\begin{array}{ll}
u_{int}
=&\ds{(\ell\varepsilon)^{1\over 2}\left[{\sqrt{ C_N^{(4)}}\over 2}e^{(N-2)\xi\over 2}
+{\sqrt{ C_N^{(4)}}\over 2}e^{(2-N)\xi\over 2}|x|^{2-N}-{4\mu  C_N^{(8)}\ell^{4-N\over N-2}\over
(N-2)^2\sqrt{ C_N^{(4)}}}e^{(N-6)\xi\over 2}
\right]}\\
&+V_{\varphi_i}+{\cal O}(\varepsilon^{1\over 2}r_\varepsilon^2)
\end{array}
\]
in $B(x_i, 2r_\e) - B(x_i, r_\e/2)$.
\end{pro}

Since we have found the solution of (\ref{eq5.2}) with the form
(\ref{eq6.4}), we have
\beq \label{eq6.36} -\Delta w_i=\lambda
(w_i+V_{\varphi_i})+(u_{p,\lambda,\xi_i}+w_i+V_{\varphi_i})^p
-u_{p,\lambda,\xi_i}^p \eeq so that in $B(x_i, r_\varepsilon)-
B(x_i, r_\varepsilon/2)$
$$
|\Delta w_i|\leq |w_i+V_{\varphi_i} |+ c
u_{p,\lambda,\xi_i}^{p-1}|w_i+V_{\varphi_i} |\leq c|w_i+V_{\varphi_i}
|
$$
Using the standard elliptic theory, we have $$
\|r_\varepsilon\nabla w_i\|_{L^\infty(B(x_i,
r_\varepsilon)- B(x_i, 3r_\varepsilon/4))} \leq c
\varepsilon^{1\over 2}(r_\varepsilon^3+ r_\varepsilon^{({2\over
p-1}+2+\delta-{2\over N+2}) {N+2\over 2}-(\delta+{2\over p-1})})
$$ Recall
$$
({2\over p-1}+2+\delta-{2\over N+2}) {N+2\over 2}>2
$$
Thus, $$ \|r_\varepsilon\nabla w_i\|_{L^\infty(B(x_i,
r_\varepsilon)- B(x_i, 3r_\varepsilon/4))} \leq c
\varepsilon^{1\over 2}\,r_\varepsilon^2
$$ By the regularity theory, for all $\alpha\in (0,1)$, \beq
\label{eq6.39} \|r_\varepsilon\,\partial_n
w_i\|_{C^{1,\alpha}(S^{N-1})} \leq c \varepsilon^{1\over
2}\,r_\varepsilon^2 \eeq
\subsection{Solutions of the nonlinear problem in
$\Omega_{ext,\varepsilon}$}

Given a $m$ functions ${\bf \phi} = (\phi_1,\ldots, \phi_m) \in
{\cal E}^m$, we now construct a family of positive solution of
(\ref{eq5.2}) in $\Omega_{ext,\varepsilon}$ which in some sense is
parameterized by ${\bf \phi}$.

\medskip

Let $\chi$ be a $C^\infty$ cut-off function defined in ${\R}^N$,
such that $\chi|_{B(0, r_0)}\equiv 1$ and  $\chi\equiv 0$ on $\R^N -
B(0,2r_0)$ and $\chi\geq 0$. Denote by $W_{\phi_i}$ the unique
harmonic extension of $\phi_i$ in $\R^N - B(x_i,r_\varepsilon)$
which decays at $\infty$. We look for a solution of (\ref{eq5.2})
in $\Omega_{ext,\varepsilon}$ of the form \beq \label{eq6.18}
u_{ext} = \sum_{i=1}^m \left(\Lambda_i \varepsilon^{1\over 2} \,
G(\cdot ,x_i) + \chi(\cdot -x_i) \, \left(W_{\phi_i} + {a_i\cdot
(\cdot -x_i)\over |\cdot -x_i|^N}\right)\right) + w_{ext} \eeq
where ${\bf a} := (a_1,\ldots ,a_m)\in (\R^N)^m$ and the function
$w_{ext}$ is assumed to be small and to satisfy
$w_{ext}|_{\partial \Omega_{ext,\varepsilon}} =0$.

\medskip
We use the maximum principle to reduce (\ref{eq5.2}) to
\beq\label{eq6.20} \left\{
\begin{array}{rllllll}
-\Delta w_{ext} & = & \lambda  \, w_{ext} + q +
Q_{\Lambda,\phi,a}(w_{ext}) &\mbox{in } &
\Omega_{ext,\varepsilon}\\[3mm]
w_{ext} & = & 0 &\mbox{on } & \partial\Omega_{ext,\varepsilon}
\end{array}
\right. \eeq where \beq \label{eq6.21} Q_{\Lambda,\phi,a}(w) =
\left| \sum_{i=1}^m \left(\Lambda_i \varepsilon^{1\over 2}G( \cdot
,x_i) +\chi( \cdot -x_i) \left(W_{\phi_i} + {a_i\cdot
(\cdot-x_i)\over |\cdot-x_i|^N}\right)\right) + w \right|^p \eeq
and where the function $q$ is given by \beq \label{eq6.22}
\begin{array}{rllll}
q(z)& = &\ds{\sum_{i=1}^m \Delta\chi(z-x_i)\left(W_{\phi_i}(z)+
{a_i\cdot (z-x_i)\over |z-x_i|^N}\right)}\\[3mm]
&+ & \ds{2  \sum_{i=1}^m \nabla\chi(z-x_i)
\cdot\nabla\left(W_{\phi_i}(z)+ {a_i\cdot (z-x_i)\over
|z-x_i|^N}\right)}\\[3mm]
& + & \ds{\lambda\sum_{i=1}^m \left(\Lambda_i \varepsilon^{1\over
2}G(z,x_i) +\chi(z-x_i)\left(W_{\phi_i}(z)+ {a_i\cdot (z-x_i)\over
|z-x_i|^N}\right)\right)}
\end{array}
\eeq Given $\Lambda_0$ and $\kappa >0$, we define
 \beq \label{eq6.23} {\cal G}=\{\Lambda\in \R^m \, : \,  |\Lambda|\leq \Lambda_0 \}
 \qquad \mbox{and} \qquad  A_{\varepsilon}=\{a\in (\R^N)^k \, : \,  |a|\leq
\varepsilon^{1\over 2}\, r_\varepsilon^N\}, \eeq Furthermore,
given $\nu \in (2-N, 3-N)$, we consider
$$
{\cal E}_{ext,\varepsilon}=\{w\in C^{0}_\nu
(\Omega_{ext,\varepsilon}) \, : \,  \| w\|_{C^{0}_\nu} \leq \kappa
\, \varepsilon^{1\over 2} \, r_\varepsilon^{2-\nu} \quad
\mbox{and} \quad w|_{\partial \Omega_{ext,\varepsilon}}=0\},
$$

For all $ {\bf a} \in A_{\varepsilon}$, ${\bf \Lambda} \in {\cal
G}$ and $\phi_j \in {\cal E}$, we estimate \beq \label{eq6.29} \|q
\|_{C^{0}_{\nu-2}(\Omega_{ext,\varepsilon})} \leq c \,
\varepsilon^{1\over 2} \, r_\varepsilon^{N} \eeq and given $w \in
{\cal E}_{ext,\varepsilon}$, we obtain with little work \beq
\label{eq6.32} \|\lambda w\|_{C^{0}_{\nu-2}}\leq c \, \|\lambda
w\|_{C^{0}_{\nu}}\leq c \, \varepsilon^{N-4\over N-2} \|
w\|_{C^{0}_{\nu}}\ \eeq and \beq \label{eq6.31}
\|Q_{\Lambda,\phi,a}(w) \|_{C^{0}_{\nu-2}} \leq  c \, (\Lambda_0^p+1+\kappa^p \varepsilon^{p-1\over 2})\,
\varepsilon^{1\over 2}  \, r_\varepsilon^{-\nu+2} \eeq Finally, we
estimate for all $w_1,w _2\in {\cal E}_{ext,\varepsilon}$ \beq
\label{eq6.33} \|Q_{\Lambda,\phi,a}(w_1) -Q_{\Lambda,\phi,a}(w_2)
\|_{C^{0}_{\nu-2}} \leq c \, \varepsilon^{2\over N+2} \, \|w_1-w_2
\|_{C^{0}_{\nu}} \eeq

The following result is standard
\begin{lem} \label{lem7.1} Assume that $\nu \in (2-N, 0)$ then for all $f\in
C^0_{\nu-2}(\Omega_{ext,\varepsilon})$, there exists $w\in
C^0_{\nu}(\Omega_{ext,\varepsilon})$ unique solution of \beq
\left\{
\begin{array}{rlll}
\Delta w & = & f & \mbox{in}~~~ \Omega_{ext,\varepsilon}\\[3mm]
       w & = & 0 & \mbox{on}~~~ \partial\Omega_{ext,\varepsilon}\,.
\end{array}
\right.
\eeq
Furthermore, there holds
$$
\|w\|_{C^0_{\nu}}\leq c\|f\|_{C^0_{\nu-2}} .
$$
\end{lem}
{\it Proof.} The existence of $w$ is straightforward and the
estimate relies on the fact that $x \rightarrow |x-x_i|^\nu$ can
be used as a barrier  in $B(x_i, r_0)- B(x_i, r_\e)$. \cqfd

\medskip

We define the map
\[
T_{\Lambda,\phi,a} : {\cal E}_{ext,\varepsilon}\longrightarrow
{\cal E}_{ext,\varepsilon}
\]
by $T_{\Lambda,\phi,a}(w) := v$ where $v$ is the solution of
\[
\Delta v =  \lambda w+q + Q_{\Lambda,\phi,a}(w).
\]

Given $\kappa >0$, it follows from the estimates (\ref{eq6.29}),
(\ref{eq6.32}) and (\ref{eq6.31}) that the mapping
$T_{\Lambda,\phi,a}$ is well defined and is a contraction,
provided $\varepsilon$ is chosen small enough, say $\epsilon \in
(0, \epsilon_0)$. In particular, this mapping has a unique fixed
point in ${\cal E}_{ext,\varepsilon}$ which yields a solution of
(\ref{eq6.20}). Therefore, we have proved the following~:
\begin{pro}
\label{pro7.8} Given ${\bf x}\in\O^m$, ${\bf a}\in (\R^N)^m$ and
${\bf\phi} \in {\cal E}^m$, there exists $u_{ext}$ positive solution of
equation (\ref{eq5.2}) in $\Omega_{ext,\varepsilon}$, satisfying
$$
u_{ext} =\phi_i +  {a_i \cdot(\cdot-x_i)\over r_\e^N}+\sum_{j=1}^m \Lambda_i \varepsilon^{1\over 2} \,
G(\cdot ,x_i)
$$
on $\partial B(x_i, r_\varepsilon)$ for all $1\leq i\leq m$ and
$u_{ext}=0$ on $\partial\O$. Furthermore, the function $u_{ext}$
can be expanded as
$$
u_{ext} = W_{\phi_i}+ {a_i \cdot(\cdot-x_i)\over r_\e^N}+\sum_{j=1}^m \Lambda_i \varepsilon^{1\over 2} \,
G(\cdot ,x_i)+{\cal O}(\e r_\e^2)
\]
in $B(x_i, 2r_\e)-B(x_i, r_\e/2)$.
\end{pro}

\medskip

Similarly, \beq \label{eq6.40}  \, \|r_\varepsilon\,\partial_n
w_{ext}\|_{C^{1,\alpha}(S^{N-1})} \leq c \varepsilon^{1\over
2}r_\varepsilon^2 \eeq where $n$ is the outside unit normal vector
on the boundary of $B(x_i, r_\varepsilon)$. In the following
consideration we will fix some $\alpha\in (0,1)$.

\subsection{The cauchy data mapping}

We explain how the free parameters in Proposition  \ref{pro7.7}
and Proposition \ref{pro7.8} can be chosen so that the functions
$u_{int,i}$ and $u_{ext}$ can be glued together to obtain a positive
solution of problem (\ref{eq5.2}) in $\Omega$.

\medskip

We set ${\bf \xi} = (\xi_1,\ldots , \xi_m)$. We want to choose the
suitable parameters
\[
\Xi : = ( {\bf x,\Lambda, \varphi,\phi,a,\xi} )
\]
so that $u_{int,i}$ and $u_{ext}$ have the same Cauchy data on
each $\partial B(x_i, r_\varepsilon)$. Once this is done, the
function defined by $u =u_{int, i}$ in $B(x_i, r_\varepsilon )$
and $u= u_{ext}$ in $\Omega_{ext,\varepsilon}$ will be ${\cal
C}^1$ and solution of (\ref{eq5.2}) away from the $\partial B(x_i,
r_{\varepsilon})$. Elliptic regularity theory will then imply that
it is a solution in $\Omega$. Moreover, it will follow from the
construction itself that $u$ has the desired behavior near each
$x_i$ and this will complete the proof of Theorem \ref{th6.1}.

\medskip

Therefore, it remains to solve, for all $i=1, \ldots ,m$, the
system \beqr \label{eq6.34} \left\{
\begin{array}{rllll}
u_{int,i} & = & u_{ext}, \\[3mm]
\label{eq6.35} \partial_n u_{int,i}  & = & \partial_n u_{ext} ,
\end{array} \right.
\eeqr on $\partial B(x_i, r_\varepsilon)$.

\medskip

We denote by $\Pi_j$ the $L^2(S^{n-1})$-projection onto
$\mbox{Span} \{e_j \}$, and
$$
\begin{array}{llll}
\Pi (\phi)  : = \ds \phi - \sum_{j=0}^N \Pi_j (\phi)
\end{array}
$$
For all $i=1, \ldots, m$, the $L^2(S^{n-1})$-projection of
(\ref{eq6.34}) over the orthogonal complement of $\mbox{Span}
\{e_0, \ldots, e_N\}$ yields the system of equations \beq
\label{eq6.3.08}
\begin{array}{rllll}
\varphi_i & = & \phi_i  + F_{i,1}(\Xi), \\[3mm]
r_\varepsilon\, \partial_n V_{\varphi_i} & = &  r_\varepsilon\,
\partial_n W_{\phi_i} + F_{i,2}(\Xi),
\end{array}
\eeq
Next, we use the expansions of Lemma \ref{lem2.5}, Corollary
\ref{cor3.1} and Corollary \ref{cor3.2} to obtain the
$L^2(S^{n-1})$-projection of (\ref{eq6.34}) over
$\mbox{Span}\{e_0\}$ \beq \label{eq6.3.09}
\begin{array}{llll}
\ds{(\ell_i C_N^{(4)} \varepsilon)^{1\over 2}\left(
{r_\varepsilon^{2-N} \, e^{-(N-2) \xi_i/2}\over 2}+ {e^{(N-2) \, \xi_i/2}\over 2}-
{4{\ell_i}^{2\over N-2} \, C_N^{(8)} \, \mu \, e^{(N-6)
\xi_i/2}\over \ell_i\, C_N^{(4)}(N-2)^2}\right)}\\[3mm]
\qquad \qquad  \qquad =\ds{\varepsilon^{1\over 2}\Lambda_i\left(
r_\varepsilon^{2-N} - H(x_i,x_i) \right) + \varepsilon^{1\over
2}\sum_{l\neq i} \Lambda_l
G(x_i,x_l)+F_{i,3}(\Xi)}, \\[3mm]
\ds{(\ell_i C_N^{(4)} \varepsilon)^{1\over 2}
\left({(2-N) \, r_\varepsilon^{2-N} \, {e^{-(N-2) \xi_i/2}
}\over 2}\right)} =\ds {\varepsilon^{1\over 2} \, \Lambda_i \, (2-N) \,
r_\varepsilon^{2-N} + F_{i,4}(\Xi)},
\end{array}
\eeq Finally, the $L^2(S^{n-1})$-projection of (\ref{eq6.34}) over
$\mbox{Span}\{e_1, \ldots, e_N\}$ yields \beq
\label{eq6.3.10}\begin{array}{rllll} \ds r_\varepsilon^{1-N} \,
{a_i } +\varepsilon^{1 \over 2} \left( r_\varepsilon \sum_{l\neq
i} \Lambda_l \, \nabla_z G(x_i,x_l) - r_\varepsilon \, \Lambda_i
 \, \nabla_z H(x_i,x_i) \right) & = & F_{i,5}(\Xi), \\[3mm]
\ds{ r_\varepsilon^{1-N} \, {a_i \, (1-N)} + \varepsilon^{1\over
2} \, \left( r_\varepsilon\sum_{l\neq i} \Lambda_l \, \nabla_z
G(x_i,x_l) - r_\varepsilon \, \Lambda_i \nabla_z
H(x_i,x_i)\right)}& = & F_{i,6}(\Xi),
\end{array}
\eeq Here $F_{i,l}(\Xi)$ for $i=1,\ldots, m$ and $l=1, \ldots ,6$
are continuous maps satisfying \beq \label{eq6.42} |F_{i,l}(\Xi)|
= {\cal O} (\varepsilon^{1\over 2} r_\varepsilon^2). \eeq

\medskip

We define "Dirichlet to Neumann map" for any
$$
{\cal S}: \Pi ({\cal C}^{2,\alpha}(S^{N-1}))
\longrightarrow
\Pi ({\cal C}^{1,\alpha}(S^{N-1}))
$$
by
\[
{\cal S} (\psi ) = r_\varepsilon \, (\partial_n V_{\psi } -
\partial_n W_{\psi}),
$$
where $V_{\psi}$ (resp. $W_{\psi }$) is the harmonic extension in
the ball $B(0,r_\e)$ (resp. in $\R^N - B(0,r_\e)$) defined in
Lemma \ref{lem6.1} and Lemma \ref{lem6.2}. It is well known that
${\cal S}$ is an isomorphism \cite{PR} the norm of whose inverse
does not depend on $\e$.

\medskip

Hence, (\ref{eq6.3.08}) (\ref{eq6.3.09}) and (\ref{eq6.3.10}) are
equivalent to the following system \beq \label{eq6.44}
\begin{array}{rllll}
\varphi_i & = & G_{i,1}(\Xi),\\[3mm]
\phi_i & = & G_{i,2}(\Xi), \\[3mm]
\xi_i & = & \ds - {2\over N-2}\log\left({2\Lambda_i \over
(\ell_i C_N^{(4)})^{1\over 2} } \right) + \varepsilon^{-{1\over 2}}
r_\varepsilon^{N-2} G_{i,3}(\Xi), \\[3mm]
a_i & = & r_\varepsilon^{N-1} G_{i,4}(\Xi), \\[3mm]
\nabla_{x_i}{\cal F}_\mu (x,\Lambda) & = & \varepsilon^{-{1\over
2}} r_\varepsilon^{-1} G_{i,5}(\Xi), \\[3mm]
\nabla_{\Lambda_i}{\cal F}_\mu (x,\Lambda) & = &
\varepsilon^{-{1\over 2}} G_{i,6}(\Xi),
\end{array}
\eeq where $G_{i,l}(\Xi)$ for all $l=1,\ldots ,6$ and for all
$i=1, \ldots ,m$ are continuous maps satisfying
$$
|G_{i,l}(\Xi)|={\cal O}(\varepsilon^{1\over 2} r_\varepsilon^2)
$$
and
$$
C_N^{(1)}={2^{4\over N-2}C_N^{(8)}\over (N-2)(C_N^{(4)})^{2\over N-2}}\qquad
C_N^{(2)}={C_N^{(4)}\over 2}.
$$
Moreover, elliptic regularity Theory shows that all $G_{i,l}(\Xi)$
are compact operators.

\medskip

Assume that $( {\bf x}^0 , \Lambda^0)$ is a non degenerate
critical point of ${\cal F}_\mu$. In particular, this implies that
$d{\cal F}_\mu$, evaluated at this point, is a local
diffeomorphism from a neighborhood of $({\bf x}^0, \Lambda^0)$ on
a neighborhood of $0$ in $\R^{m(N+1)}$. Using this we can write
formally the system (\ref{eq6.44}) as
\[
\Xi = \Phi(\Xi),
\]

We set $\xi_i^0 : = - {2\over N-2}\log\left({2\Lambda_i^0 \over
(\ell_i C_N^{(4)})^{1\over 2}}\right)$ for all $i=1,\ldots ,m$. We
consider the set
$$
{\cal A}=\overline{B(({\bf x}^0 , \Lambda^0 ),\varepsilon_1)} \times
{\cal E}^{2k}\times A_\varepsilon\times \overline{B({\bf \xi}^0
,\varepsilon_1)}
$$
where $\varepsilon_1$ is some fixed small positive number. It follows from the above analysis that $\Phi:{\cal A}\rightarrow
{\cal A}$ is a continuous compact map. According to Schauder fixed
point theorem, $\Phi$ has a fixed point in ${\cal A}$. This
completes the proof of Theorem 1.

\medskip

\begin{rem}
If $N=4$, we have the similar results. In this case, we take
$\lambda=\mu \frac{1}{\log{1/\varepsilon}}$.
\end{rem}

%{\small \noindent

%Y. Ge ( ge@univ-paris12.fr )

%\noindent
%Laboratoire  d'Analyse et de Math\'ematiques Appliqu\'ees, 
%CNRS UMR 8050, 
%D\'epartement de Math\'ematiques, 
%Universit\'e Paris XII-Val de Marne, 61 avenue du G\'en\'eral de Gaulle, 
%94010 Cr\'eteil Cedex, France}
%
%
%
%
%\bigskip
%
%{\small \noindent
%
%R. Jing (jing@univ-paris12.fr)
%
%
%\noindent
%Laboratoire  d'Analyse et de Math\'ematiques Appliqu\'ees, 
%CNRS UMR 8050, 
%D\'epartement de Math\'ematiques, 
%Universit\'e Paris XII-Val de Marne, 61 avenue du G\'en\'eral de Gaulle, 
%94010 Cr\'eteil Cedex, France\\
%and\\
%\noindent
%Department of Mathematics, 
%East China Normal University, 
%200062 Shanghai, P.R. of China}
%
%\bigskip
%
%{\small \noindent
%
%F. Pacard ( pacard@univ-paris12.fr )
%
%\noindent
%Laboratoire  d'Analyse et de Math\'ematiques Appliqu\'ees, 
%CNRS UMR 8050, 
%D\'epartement de Math\'ematiques, 
%Universit\'e Paris XII-Val de Marne, 61 avenue du G\'en\'eral de Gaulle, 
%94010 Cr\'eteil Cedex, France}

 \end{document}